\documentclass[14pt]{extarticle}

\usepackage[margin=1in]{geometry}  % set the margins to 1in on all sides
\usepackage{graphicx}              % to include figures
\usepackage{amsmath} \allowdisplaybreaks   % great math stuff
\usepackage{amsfonts}              % for blackboard bold, etc
\usepackage{amsthm}                % better theorem environments
\usepackage{CJK}
\usepackage{float}
\usepackage{amssymb}
\usepackage{tikz}
\usepackage{setspace}
\usepackage{indentfirst}
\usepackage{footmisc} 
\usepackage{booktabs}
\usepackage{array}
\usepackage{float}
\usepackage[lined,linesnumbered]{algorithm2e}
\usepackage{bm}

\SetAlFnt{\normalsize}
\SetAlCapFnt{\normalsize\bfseries}
\SetAlCapNameFnt{\normalsize\bfseries}
\SetInd{0.7em}{1.2em}

\SetCommentSty{myCommentStyle}

\SetKwComment{tcp}{\(\triangleright\)\ }{}
\SetKwComment{tcc}{\(\triangleright\)\ }{}

\makeatletter
\let\oldnl\nl
\newcommand{\nonl}{\renewcommand{\nl}{\let\nl\oldnl}}
\makeatother

\usepackage[backend=biber, 
            style=authoryear,
            maxcitenames=1,
            mincitenames=1, 
            sorting=nyt, 
            dashed=false]
{biblatex}   % 按姓名、年份、标题排序

\DeclareFieldFormat{citehyperref}{\textmd{#1}}
\DeclareFieldFormat{labeldesignator}{\textmd{#1}}
\DeclareFieldFormat{author}{\textmd{#1}}
\DeclareFieldFormat{year}{\textmd{#1}}

\AtEveryCite{%
  \let\bfseries\mdseries
}

\DeclareFieldFormat[article,book,inbook,incollection,thesis]{title}{#1}

\DeclareFieldFormat{journal}{#1.} % 期刊名称后加点号
\renewbibmacro{in:}{} 

\renewbibmacro*{volume+number+eid}{%
  \printfield{volume}% 打印卷号
  \iffieldundef{number}
    {}
    {\printtext[parens]{\printfield{number}}}% 用括号包裹期号
  \printtext{:}% 在卷号+期号后添加冒号
  \setunit{\addspace}% 添加空格
}

\DeclareFieldFormat{pages}{#1} % 页码直接打印，不加其他符号

\SetKwInOut{INPUT}{input}
\SetKwInOut{OUTPUT}{output}
\SetKw{RETURN}{return}

\newtheorem{lem}{Lemma}

\newtheorem{conj}{Conjecture}

\newtheorem{defn}{Definition}
\theoremstyle{remark}

\newtheorem{example}{Example}

\providecommand{\keywords}[1]
{
  \small	
  \textbf{\textit{Keywords---}} #1
}

\usepackage{xcolor}
\usetikzlibrary{arrows.meta,calc}

\definecolor{panelblue}{RGB}{241,247,253}
\definecolor{arrowblue}{RGB}{55,113,195}
\definecolor{rootpurple}{RGB}{205,204,255}
\definecolor{activeblue}{RGB}{207,218,255}
\definecolor{processedgray}{RGB}{215,215,215}
\definecolor{encodingred}{RGB}{235,35,35}
\definecolor{graytext}{RGB}{145,145,145}
\definecolor{navy}{RGB}{20,48,88}

\title{A complete classification for  n-species Lotka-Volterra stably dissipative systems based on a generalized AHU algorithm}

\author{
  Yajia Du\footnotemark[1] 
  \and
  Zhengyi Lu\footnotemark[2] 
  \and
  Jinya Xu\footnotemark[3]
  \and
  Ming Yang\footnotemark[2] 
}

\date{} % Comment this line to show today's date

\begin{document}
\onehalfspacing
\renewcommand{\thefootnote}{\fnsymbol{footnote}}  % 使用符号作为脚注标记
\maketitle
\begingroup
\setlength{\parindent}{0pt}
\renewcommand{\footnotelayout}{\raggedright}
%\footnotetext[1]{No.9 Middle School of Nanchong, No.11 Fengdeng Road, Shunqing District, Nanchong 637503, Sichuan, China} 
%\footnotetext[2]{School of Mathematics and Software Science, Sichuan Normal University, Chengdu 610101, Sichuan, China}
%\footnotetext[3]{School of Intelligent Science and Engineering, Chengdu Neusoft University, Dujiangyan Chengdu 611844, Sichuan, China. Email: xujinya@nsu.edu.cn. Corresponding author.}
\footnotetext[1]{No. 9 Middle School of Nanchong, No. 11 Fengdeng Road, Shunqing District, Nanchong 637503, Sichuan, China.}

\footnotetext[2]{School of Mathematics and Software Science, Sichuan Normal University, Chengdu 610101, Sichuan, China.}

\footnotetext[3]{School of Intelligent Science and Engineering, Chengdu Neusoft University, Dujiangyan Chengdu 611844, Sichuan, China. Email: xujinya@nsu.edu.cn. Corresponding author.}
\endgroup

\abstract{
This paper proposes a generalized AHU algorithm based on tree isomorphism theory and applies it to the complete classification of $n$-dimensional stably dissipative Lotka-Volterra systems.}
\bigskip\bigskip

\bigskip\bigskip

\section{Introduction}

\bigskip\bigskip

Consider the following $n$-species Lotka-Volterra system
\begin{equation}
\label{1}
\dot{x}_i=x_i(r_i+\sum_{j=1}^n a_{ij}x_j), \quad i=1,...,n.
\end{equation}
Here, $r_i$ and $a_{ij}$ are real constants. We focus on solutions residing in $\mathbb{R}_+^n$. It is well established that the global dynamics of system (1) is governed by its interaction matrix $A = (a_{ij})_{n \times n}$ [\cite{Duarte1998, Goh1977,Hsu1978,Lu1995,Redheffer1985,Redheffer1989,Zhao2010}].
%[\cite{Duarte1998}, \cite{Goh1977}, \cite{Hsu1978}, \cite{Lu1995}, \cite{Redheffer1985}, \cite{Redheffer1989}, \cite{Zhao2010}].  
For matrix $A$, there is an important concept which determines the dynamical behavior of the system[\cite{Cross1978}, \cite{Goh1977}, \cite{Hsu1978}, \cite{Redheffer1981}].

\begin{defn}
Matrix $A$ is Volterra-Lyapunov $($semi-$)$ stable, if there is a positive diagonal matrix $C=$ \textnormal{diag}$(c_1,...,c_n)$
such that $CA+A^TC$ is negative $($semi-$)$ definite.
\end{defn}

Suppose that system(1) admits a unique positive equilibrium $x^*=(x^*_1,...,x^*_n)$. Then system(1) can be reformulated  as

\begin{equation}
\dot{x}_i=x_i \sum_{j=1}^n a_{ij}(x_j-x^*_j), \quad i=1,...,n.
\end{equation}

By employing the Volterra-Lyapunov function

\begin{equation}
V(x)= \sum_{i=1}^n c_{i}(x_i-x^*_i-x^*\ln \frac{x_i}{x^*_i}),
\end{equation}
and applying LaSalle's invariance principle[\cite{LaSalle1960,Harrison1979,Hsu1978,Krikorian1979}], we obtain the following results:

\begin{lem}
If $A$ is Volterra-Lyapunov semi-stable, then every solution of system \textnormal{(2)} converges to
 the maximal invariant set $M$ contained in the set
\begin{equation}
E=\{x\in {\rm int \mathbb{R}^n_+} | (x-x^*)^T(CA+A^TC)(x-x^*)=0 \}.
\end{equation}
\end{lem}
\noindent Here, $M$ denotes the LaSalle's invariant set of system (2). Specifically, if the matrix $A$ is Volterra-Lyapunov stable, $M$ reduces to a singleton $\{x^*\}$, thereby implying that $x^*$ is globally asymptotically stable.

However, when $A$ is merely Volterra-Lyapunov semi-stable (but not stable), the LaSalle's invariant set $M$ contained in $E$ typically extends well beyond a single point $x^*$. In such cases, characterizing the geometric or algebraic structure of $M$ becomes a formidable challenge.

To address this, several intermediate notions of stability-positioned between Volterra-Lyapunov semi-stability and stability have been explored in the literature. These stability conditions are strictly weaker than Volterra-Lyapunov stability yet stronger than semi-stability.

In this context, \cite{Lu1995} proposed a conjecture regarding these properties, which was subsequently proven to hold for the case $n=3$.

\begin{conj}
If $A$ is Volterra-Lyapunov semi-stable and diag$(x^*)A$ is stable, then $x^*$ is globally stable for
system \textnormal{(2)}.
\end{conj}

Following the foundational work of Samuelson[\cite{Samuelson1947}], Lancaster [\cite{Lancaster1962}, \cite{Lancaster1964}], Gorman [\cite{Gorman1964}], Quirk and Ruppert [\cite{Quirk1965}] provided a detailed qualitative stability analysis for general matrices. Subsequently, Jeffries [\cite{Jeffries1988}, \cite{Jeffries1974}] established a complete set of criteria for qualitative stability using graph methods and color test notation.
Based on a theorem by Redheffer[\cite{Redheffer1989}], Lu and Takeuchi[\cite{Lu1994}] proved that the
qualitative stability of the interaction matrix $A$ of system (2) implies its global stability.
Furthermore, any irreducible matrix possessing qualitative stability must be Volterra-Lyapunov semi-stable.

In a series papers [\cite{Liu1991}, \cite{Lu1988}, \cite{Lu1989}], Lotka-Volterra prey-predator chain systems were extensively investigated.

\begin{defn}
 Lotka-Volterra system \textnormal{(2)} is defined as a prey-predator cycle, if its interaction matrix $A = (a_{ij})$ satisfies the following conditions:
 
\textnormal{1)}  \quad $a_{ii}\leq 0$ for all $i = 1, \dots, n$, with at least one strict inequality holding;

\textnormal{2)}  \quad $a_{ij}a_{ji} \leq 0$ and $a_{ij}=0$ if and only if $a_{ji}=0$, $i \not=j$, $i,j=1,…,n$;

\textnormal{3)}  \quad there are distinct indices $i_1, \cdots,i_m (m\geq 3)$ such that none of the element $a_{i_1i_2},a_{i_2i_3}, \cdots, a_{i_mi_1}$
vanish.
\end{defn}
A system that does not satisfy the requirements of a prey-predator cycle is referred to as a prey-predator chain system.

It can be verified that the interaction matrix of a prey-predator chain system is necessarily Volterra-Lyapunov semi-stable.In particular, Lotka-Volterra prey-predator chain systems are globally stable for dimensions $n(\leq 4)$ [\cite{Redheffer1982}]. Regarding these systems, Chen [\cite{Lu1989}] proposed the following conjecture:

\begin{conj}
Every locally asymptotically stable equilibrium of the n-dimensional Lotka-Volterra 
prey-predator chain system \textnormal{(2)} is globally stable.
\end{conj}

It was previously demonstrated that the Chen conjecture holds true for the n-species food chain system (2) where $n \le 8$[\cite{Yang2019}].

A concept slightly broader than qualitative stability is stable dissipation, introduced by Redheffer[\cite{Redheffer1982}]. This concept is based on the notion that a perturbation of a matrix $A=(a_{ij})_{n\times n}$ is defined as any matrix  $\bar{A}=(\bar{a}_{ij})_{n\times n}$ such that $\bar{a}_{ij}=0$ if and only if $a_{ij}=0$.

\begin{defn}
(Stably Dissipative Matrix)  The interaction matrix $A$ of System \textnormal{(2)} is called stably dissipative,  denoted by $A\in \mathbb{S_D}$, if every sufficiently small perturbation $\bar{A}$ of $A$
 is Volterra-Lyapunov semi-stable: that is, there exists $\epsilon >0$ such that $\max_{ij}|a_{ij}-\bar{a}_{ij}|<\epsilon$
$\Rightarrow \bar{A}$ is Volterra-Lyapunov semi-stable.
\end{defn}

The perturbations considered here are required to preserve the zero pattern of \(A\), that is, \(\bar{a}_{ij}=0\) if and only if \(a_{ij}=0\). A Lotka-Volterra system with a stably dissipative interaction matrix is called a stably dissipative system.

 Section 3 reviews the properties of stably dissipative matrix and and presents relevant known results, together with several key theorems and lemmas. Section 3 revisits the classical AHU algorithm for encoding rooted and unrooted trees and develops a generalized AHU algorithm that extends the encoding framework to unrooted colored trees with black and white nodes.  Build on this generalization, Section 4 presents a comprehensive classification algorithm for n-dimensional stably dissipative Lotka-Volterra systems and reports computational results on tree isomorphism, with the computations extending up to $n=15$ for unrooted colored trees and $n=21$ for uncolored trees. Finally, concluding remarks are offered in Section 5.

\section{Criteria for stably dissipative matrice}

 If $A \in \mathbb{S_D}$, then $A$ is Volterra-Lyapunov semi-stability, i.e.,  there is a positive diagonal matrix $C =$diag$(c_1, \dots, c_n)$ such that $C A + A^T C$ is semi-negative definite. Furthermore,  if[\cite{Redheffer1982}]
  \begin{equation}
  \label{2.1}
    (x - x^*) (C A + A^T C) (x - x^*)^T = 0
  \end{equation}
then
\begin{equation}
    \label{2.2}
    a_{ii} (x - x^*)_i = 0,\ i = 1, \dots, n.
  \end{equation}
\bigskip
Here, (\ref{2.1}) and (\ref{2.2}) imply that in the LaSalle's invariant set when $a_{ii}\not= 0$, $x_i=x_i^*$.   If $a_{ii}<0$ ($i=1,\ldots,n$), then $A\in \mathbb{S_D}$ if and only if $A$ is Volterra-Lyapunov stable[\cite{Redheffer1982}, \cite{Zhao2010}].
\bigskip

Clearly, in this special case for  $a_{ii}<0$ ($i=1,\ldots,n$), the dynamical behavior for system (2) is simple, i.e., globally stable.  Therefore the interesting phenomena will appear when some of $a_{ii}$ being zero. To deal with this situation, we define as [~\cite{Redheffer1982}]
 the interaction graph $G(A)$ of the system (1) or (2) with the $n \times n$ interaction matrix $A = (a_{i,j})_{n \times n}$ as follows:

   $G(A)$ consists of $n$ nodes representing $n$ populations together with edges representing the interrelations between populations.

   A black dot $\bullet$ (when $a_{i,i}<0$) or white dot $\circ$ (when $a_{i,i}=0$) will be marked on each node.

   There is an edge between $i$ and $j$  whenever $a_{i,j} \neq 0$ and $a_{j,i} \neq 0$.

In[~\cite{Hofbauer1998}], it is shown that if $G(A)$ is free of loop, then $A$ is Volterra-Lyapunov stable if and only if $-A$ is a P matrix.  In [~\cite{Redheffer1982}], it is shown that if $A\in \mathbb{S_D}$, $i$ and $j$ are the nodes connected to $G(A)$, then $a_{ii}a_{jj}>a_{ij}a_{ji}$.

The following strong (weak) connection concept is important to classify the graph of the interaction matrix.

\begin{defn}
\textnormal{[\cite{Redheffer1982}]}\quad A connection $[i,j]$ in $G(A)$ is called a strong (weak) connection, if $a_{ij}a_{ji}<0$ \textnormal{(}$a_{ij}a_{ji}=0$\textnormal{)}.
\end{defn}

It is shown [\cite{Zhao2010}] that if $A\in \mathbb{S_D}$, then each loop in $G(A)$ has at least one strong connection.  Furthemore, we have

\begin{lem}
\textnormal{[\cite{Zhao2010}]}\quad If the graph $G(A)$ corresponding to $A$ has only one black node $\bullet$, then $A \in \mathbb{S_D}$ if and only if

\textnormal{(i)}\ $G(A)$ is free of loop and

\textnormal{(ii)}\ If $a_{ij}\neq 0$, then $a_{ij}a_{ji}<0$ is true for all $i\neq j$.
\end{lem}

\begin{lem}\textnormal{[\cite{Zhao2010}]}
\quad Suppose $G(A)$ has no less than two black nodes $\bullet$, i.e., $a_{ii}<0$ $(i=1,\ldots,k$; $k\geq 2)$ and $a_{jj }=0$ $(j=k+1,\ldots, n)$. Let $M=(a_{ij})_{k\times k}$ $(1\leq i,\ j\leq k)$ and  $\bar G(A)$ denote the graphs $G(A)$ in which all the links $[i,j]$ are removed $(i\leq k$, $j\leq k)$ (actually, $\bar G(A)$ is obtained by removing all strongly connected graphs from $G(A))$, then $A\in \mathbb{S_D}$ if and only if.

\textnormal{(i)}\ $\bar G(A)$ is free of loop;

\textnormal{(ii)}\ There is a positive diagonal array $C=$\textnormal{diag}$(c_1,\ldots,c_n)$ such that $C_0M+M^TC_0$ is negative and $c_ia_{ij}+c_ja_{ji}=0$ $(i>k$, $j>k)$, here $C_0=$\textnormal{diag}$(c_1,\ldots,c_k)$.
\end{lem}

According to the two lemmas above, the matrix $A$ is stably dissipative if and only if the graph $\bar{G}(A)$---formed by removing all strongly connected components (or strong edges) from $G(A)$---is a tree[\cite{Zhao2010}].  The sole exception occurs when all diagonal elements of $A$ are negative, which implies that all edges in $G(A)$ are strongly connected. Consequently, the problem of system classification based on the interaction matrix $A$ can be mapped onto the topological classification of its associated graph $G(A)$. To this end, we have:

\begin{defn}
The graph $G(A)$ of matrix $A$ is called the maximal stably dissipative graph, if $A \in \mathbb{S_D}$ but $\bar A \notin \mathbb{S_D}$. 
Here $\bar A$ is obtained by adding any links between the nodes of G(A).
\end{defn}

Based on this definition, once the set of all maximal stably dissipative graphs for $n$-order matrices is determined, the stably dissipative graph $G(A)$ of any $n$-order tolerance matrix $A$ can be retrieved by removing $k$ ($k \geq 0$) strong connections from a maximal graph. Therefore, classifying maximal stably dissipative graphs is essentially equivalent to classifying all stably dissipative graphs. The computational process is straightforward: all stably dissipative graphs can be exhausted by removing subsets of strong connections from the identified maximal graphs.

A complete classification algorithm for ~$n$-dimensional ~Lotka-Volterra system with maximal stably dissipative matrix is proposed in section 3. Theoretically, for a given ~$n$, the classification of the maximal stably dissipative graph for a $n \times n$ matrix will be determined by a computer program based on the algorithm. A result up to $ n=15$ is presented. In fact,  the classification result of the maximal stably dissipative graph is the same as the known conclusion for ~$n=4$ and $5$[\cite{Zhao2010}].  For ~$n=6$,  [\cite{Yu2011}] obtained  ~73 classes maximal stably dissipative graphs, but in our case, it should be ~76 classes. For  ~$n=7$,  [\cite{Zhang2011}]  found ~229 classes maximal stably dissipative graphs, in our case, it is ~235 classes. 
\section{Canonical Encoding of Bichromatic Trees}
This section presents an algorithm for identifying maximally stably dissipative graphs. We extend the AHU algorithm to accommodate graphs whose nodes are assigned binary colors, namely, black and white. We begin by introducing several fundamental concepts.
\bigskip

\begin{defn}\textnormal{[\cite{Valiente2002}]}
The number of edges incident to a node in $G(A)$ is called the degree of that node, A node with degree 1 is called a leaf.
\end{defn}

We assign each node of a tree one of two colors, black or white. A fundamental property of trees is that removing a leaf from an $n$-node tree yields a tree with $n-1$ nodes. Conversely, any $n$-node tree can be constructed by attaching a leaf to an $(n-1)$-node tree. This inductive property provides a systematic framework for constructing all stably dissipative trees.

Apart from the inclusion of node coloring (bichromatic node), these trees retain the same structural properties as conventional trees in graph theory. Accordingly, graph isomorphism must account for both the topological structure and the color of each node. Although the AHU algorithm works well for uncolored trees, its direct application  to bichromatic trees presents considerable challenges due to the substantially increased complexity of distinguishing non-isomorphic structures. For instance, there are only three non-isomorphic uncolored trees with five nodes; however, when node coloring is considered, the number of distinct isomorphism classes increases to 26.

In the following section, we develop a generalization of the AHU algorithm for determining isomorphism among bichromatic trees. We begin by formally introducing the distinctions among rooted trees, unrooted trees and colored trees.

\subsection{The Generalized AHU Algorithm}

The original AHU algorithm employs a recursive nested-parenthesis encoding scheme to generate canonical strings for trees. In this study, we propose a generalized AHU algorithm, with the initial focus on a unique canonical encoding for rooted trees, as defined in Section 3.2 (Definitions 7 and 8). For unrooted trees, we use their graph centers to determine the root structure, thereby reducing the problem to rooted-tree encoding. It can be formally proved that the proposed encoding is canonical: two trees are isomorphic if and only if they have the same encoding.

Note that the center of a tree consists of either a single node or a pair of adjacent nodes. To address the latter case, we generate separate encodings for each node in the central pair and concatenate them in a lexicographical (sorted) order. This approach maintains the uniqueness of the encoding for all unrooted trees. Finally, we incorporate node coloring by assigning positive and negative signs to represent black and white nodes, respectively. This extension ensures that the final encoding remains a canonical form for bicolored trees.

\subsection{Encoding of Rooted Trees}
\begin{defn}
\label{def01}
Let $T=(V,E,r)$ be a rooted tree with $n$ nodes, where $V$ is the set of nodes, $E$ is the set of edges, and $r$ is the root of the tree. For any node $v \in V$, its encoding is defined as the following sequence of integers,
\begin{eqnarray}
     \label{08}
    code(v)=\left\{
    \begin{array}{ll}
        [size(v),op(code(w_1)),\cdots,op(code(w _k))] &  \mbox{v is a non-leaf node}, \\
        \left[ 1 \right] & \mbox{v is a leaf node}.
    \end{array}
    \right.
\end{eqnarray}
where $size(v)$ denotes the total number of nodes in the subtree rooted at $v$ in $T$, and $w_1, \cdots, w_k$ are the children of $v$, arranged in ascending order according to their encodings. Here, $op(code(w_j)) (j=1, \cdots, k)$, denotes the sequence of integers obtained from the encoding of $w_j$ by removing the outermost pair of square brackets. The function \(\operatorname{op}(\cdot)\) is the op function in Maple, which is used to extract the operands of an expression.
\end{defn}
\smallskip
\begin{defn}
Let $T=(V,E,r)$ be a rooted tree with $n$ nodes, where $V$ is the set of nodes, $E$ is the set of edges, and $r$ is the root of the tree. The encoding of tree $T$ is defined as the encoding of its root $r$, that is,
\begin{equation}
code(T) = code(r).
\end{equation}
\end{defn}

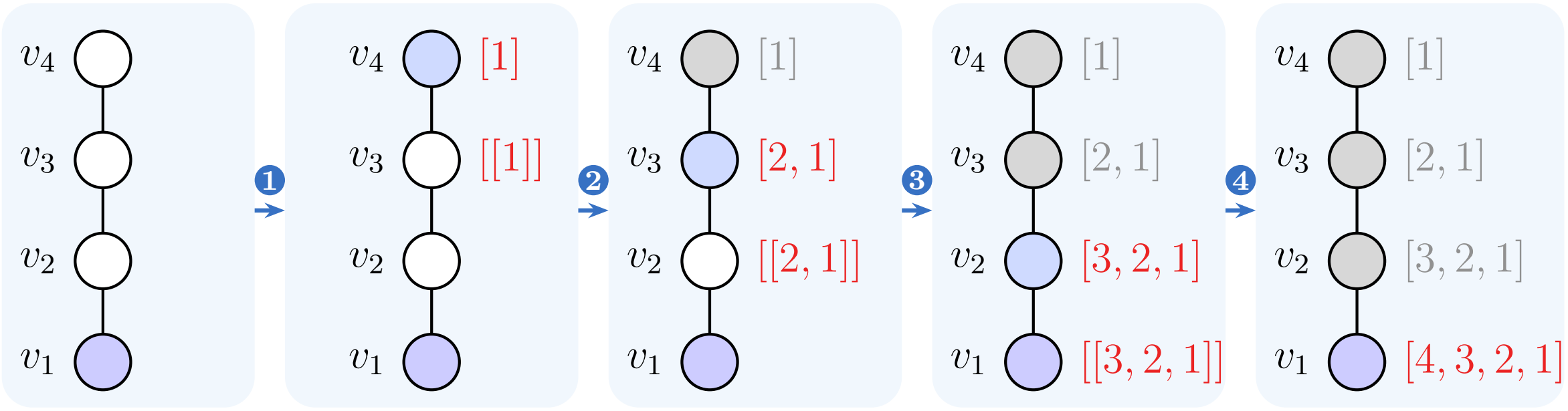
\begin{figure}[htbp]
\centering
% \resizebox{\textwidth}{!}{%
\begin{tikzpicture}[
    font=\normalsize,
    vertex/.style={
        circle,
        draw=black,
        line width=0.8pt,
        minimum size=5.5mm,
        inner sep=0pt
    },
    treeedge/.style={
        draw=black,
        line width=0.8pt
    },
    steparrow/.style={
        -{Stealth[length=2mm,width=1.5mm]},
        draw=arrowblue,
        line width=1.2pt
    },
    stepnumber/.style={
        circle,
        fill=arrowblue,
        text=white,
        font=\bfseries\scriptsize,
        minimum size=3mm,
        inner sep=0pt
    },
    panel/.style={
        rounded corners=3mm,
        fill=panelblue,
        draw=none
    },
    encoding/.style={
        font=\small,
        text=encodingred,
        anchor=west
    },
    oldencoding/.style={
        font=\small,
        text=graytext,
        anchor=west
    },
    vlabel/.style={
        font=\small,
        anchor=east
    }
]

% Horizontal positions of the five panels
\def\xA{1.5}
\def\xB{4.5}
\def\xC{7.7}
\def\xD{10.9}
\def\xE{14.1}

% Vertical positions of the four vertices
\def\yFour{4.0}
\def\yThree{3.0}
\def\yTwo{2.0}
\def\yOne{1.0}

\def\shift{0.45}

% --------------------------------------------------
% Panel (a): Original tree
% --------------------------------------------------
\fill[panel]
    (\xA-1.25,0.55) rectangle (\xA+1.25,4.55);

\draw[treeedge]
    (\xA-\shift+0.2,\yFour) -- (\xA-\shift+0.2,\yThree)
    (\xA-\shift+0.2,\yThree) -- (\xA-\shift+0.2,\yTwo)
    (\xA-\shift+0.2,\yTwo) -- (\xA-\shift+0.2,\yOne);

\node[vertex,fill=white] at (\xA-\shift+0.2,\yFour) {};
\node[vertex,fill=white] at (\xA-\shift+0.2,\yThree) {};
\node[vertex,fill=white] at (\xA-\shift+0.2,\yTwo) {};
\node[vertex,fill=rootpurple] at (\xA-\shift+0.2,\yOne) {};

\node[vlabel] at (\xA-\shift+0.2-0.28,\yFour) {$v_4$};
\node[vlabel] at (\xA-\shift+0.2-0.28,\yThree) {$v_3$};
\node[vlabel] at (\xA-\shift+0.2-0.28,\yTwo) {$v_2$};
\node[vlabel] at (\xA-\shift+0.2-0.28,\yOne) {$v_1$};

% --------------------------------------------------
% Panel (b): After processing v_4
% --------------------------------------------------
\fill[panel]
    (\xB-1.45,0.55) rectangle (\xB+1.45,4.55);

\draw[treeedge]
    (\xB,\yFour) -- (\xB,\yThree)
    (\xB,\yThree) -- (\xB,\yTwo)
    (\xB,\yTwo) -- (\xB,\yOne);

\node[vertex,fill=activeblue] at (\xB,\yFour) {};
\node[vertex,fill=white] at (\xB,\yThree) {};
\node[vertex,fill=white] at (\xB,\yTwo) {};
\node[vertex,fill=rootpurple] at (\xB,\yOne) {};

\node[vlabel] at (\xB-0.28,\yFour) {$v_4$};
\node[vlabel] at (\xB-0.28,\yThree) {$v_3$};
\node[vlabel] at (\xB-0.28,\yTwo) {$v_2$};
\node[vlabel] at (\xB-0.28,\yOne) {$v_1$};

\node[encoding] at (\xB+0.28,\yFour) {$[1]$};
\node[encoding] at (\xB+0.28,\yThree) {$[[1]]$};

% --------------------------------------------------
% Panel (c): After processing v_3
% --------------------------------------------------
\fill[panel]
    (\xC-1.45,0.55) rectangle (\xC+1.45,4.55);

\draw[treeedge]
    (\xC-\shift,\yFour) -- (\xC-\shift,\yThree)
    (\xC-\shift,\yThree) -- (\xC-\shift,\yTwo)
    (\xC-\shift,\yTwo) -- (\xC-\shift,\yOne);

\node[vertex,fill=processedgray] at (\xC-\shift,\yFour) {};
\node[vertex,fill=activeblue] at (\xC-\shift,\yThree) {};
\node[vertex,fill=white] at (\xC-\shift,\yTwo) {};
\node[vertex,fill=rootpurple] at (\xC-\shift,\yOne) {};

\node[vlabel] at (\xC-\shift-0.28,\yFour) {$v_4$};
\node[vlabel] at (\xC-\shift-0.28,\yThree) {$v_3$};
\node[vlabel] at (\xC-\shift-0.28,\yTwo) {$v_2$};
\node[vlabel] at (\xC-\shift-0.28,\yOne) {$v_1$};

\node[oldencoding] at (\xC-\shift+0.28,\yFour) {$[1]$};
\node[encoding] at (\xC-\shift+0.28,\yThree) {$[2,1]$};
\node[encoding] at (\xC-\shift+0.28,\yTwo) {$[[2,1]]$};

% --------------------------------------------------
% Panel (d): After processing v_2
% --------------------------------------------------
\fill[panel]
    (\xD-1.45,0.55) rectangle (\xD+1.45,4.55);

\draw[treeedge]
    (\xD-\shift,\yFour) -- (\xD-\shift,\yThree)
    (\xD-\shift,\yThree) -- (\xD-\shift,\yTwo)
    (\xD-\shift,\yTwo) -- (\xD-\shift,\yOne);

\node[vertex,fill=processedgray] at (\xD-\shift,\yFour) {};
\node[vertex,fill=processedgray] at (\xD-\shift,\yThree) {};
\node[vertex,fill=activeblue] at (\xD-\shift,\yTwo) {};
\node[vertex,fill=rootpurple] at (\xD-\shift,\yOne) {};

\node[vlabel] at (\xD-\shift-0.28,\yFour) {$v_4$};
\node[vlabel] at (\xD-\shift-0.28,\yThree) {$v_3$};
\node[vlabel] at (\xD-\shift-0.28,\yTwo) {$v_2$};
\node[vlabel] at (\xD-\shift-0.28,\yOne) {$v_1$};

\node[oldencoding] at (\xD-\shift+0.28,\yFour) {$[1]$};
\node[oldencoding] at (\xD-\shift+0.28,\yThree) {$[2,1]$};
\node[encoding] at (\xD-\shift+0.28,\yTwo) {$[3,2,1]$};
\node[encoding] at (\xD-\shift+0.28,\yOne) {$[[3,2,1]]$};

% --------------------------------------------------
% Panel (e): Final encoding
% --------------------------------------------------
\fill[panel]
    (\xE-1.45,0.55) rectangle (\xE+1.6,4.55);

\draw[treeedge]
    (\xE-\shift,\yFour) -- (\xE-\shift,\yThree)
    (\xE-\shift,\yThree) -- (\xE-\shift,\yTwo)
    (\xE-\shift,\yTwo) -- (\xE-\shift,\yOne);

\node[vertex,fill=processedgray] at (\xE-\shift,\yFour) {};
\node[vertex,fill=processedgray] at (\xE-\shift,\yThree) {};
\node[vertex,fill=processedgray] at (\xE-\shift,\yTwo) {};
\node[vertex,fill=rootpurple] at (\xE-\shift,\yOne) {};

\node[vlabel] at (\xE-\shift-0.28,\yFour) {$v_4$};
\node[vlabel] at (\xE-\shift-0.28,\yThree) {$v_3$};
\node[vlabel] at (\xE-\shift-0.28,\yTwo) {$v_2$};
\node[vlabel] at (\xE-\shift-0.28,\yOne) {$v_1$};

\node[oldencoding] at (\xE-\shift+0.28,\yFour) {$[1]$};
\node[oldencoding] at (\xE-\shift+0.28,\yThree) {$[2,1]$};
\node[oldencoding] at (\xE-\shift+0.28,\yTwo) {$[3,2,1]$};
\node[encoding] at (\xE-\shift+0.28,\yOne) {$[4,3,2,1]$};

% --------------------------------------------------
% Arrows and step numbers
% --------------------------------------------------
\draw[steparrow] (2.75,2.5) -- (3.05,2.5);
\node[stepnumber] at (2.90,2.8) {1};

\draw[steparrow] (5.95,2.5) -- (6.25,2.5);
\node[stepnumber] at (6.10,2.8) {2};

\draw[steparrow] (9.15,2.5) -- (9.45,2.5);
\node[stepnumber] at (9.30,2.8) {3};

\draw[steparrow] (12.35,2.5) -- (12.65,2.5);
\node[stepnumber] at (12.50,2.8) {4};

% --------------------------------------------------
% Subcaptions
% --------------------------------------------------
\node[
    font=\scriptsize,
    text=navy,
    align=center,
    text width=3.5cm
] at (\xA,0.05)
    {(a) Original tree};

\node[
    font=\scriptsize,
    text=navy,
    align=center,
    text width=3.5cm
] at (\xB-0.2,0.05)
    {(b) Processing $v_4$};

\node[
    font=\scriptsize,
    text=navy,
    align=center,
    text width=2.7cm
] at (\xC-0.2,0.05)
    {(c) Processing $v_3$};

\node[
    font=\scriptsize,
    text=navy,
    align=center,
    text width=2.7cm
] at (\xD-0.2,0.05)
    {(d) Processing $v_2$};

\node[
    font=\scriptsize,
    text=navy,
    align=center,
    text width=2.7cm
] at (\xE-0.2,0.05)
    {(e) Final encoding};

\end{tikzpicture}
% }

\caption{Step-by-step encoding of a rooted tree with four nodes.}
\label{fig-01}

\end{figure}

\begin{example}
Consider the rooted tree with four nodes and root $v_1$, as shown in Figure \ref{fig-01}.  The encoding process starts at the leaf node $v_4$ and proceeds upward toward the root. Although $v_1$ has degree one, it is not treated as a leaf in this encoding process, because it serves as the root. The encoding steps are as follows:
\end{example}
\begin{itemize}
    \item The encoding of $v_4$ is obtained as [1], where 1 represents the total number of nodes in the subtree rooted at $v_4$. This encoding is then passed to its parent node $v_3$ as [[1]] (temporarily stored, not the final encoding);
    \item Considering the removal of the subtree rooted at the previous leaf node $v_4$, we now select a new leaf node $v_3$. The total number of nodes in the subtree rooted at $v_3$ is 2. We then concatenate this with the encoding [1] of the previously temporarily stored child node $v_4$, resulting in the encoding [2,1] for node $v_3$. Simultaneously, we pass this encoding to its parent node $v_2$, becoming [[2,1]];
    \item Continuing as above, we now select another new leaf node $v_2$. The total number of nodes in the subtree rooted at $v_2$ is 3. We then concatenate this with the encoding [2,1] of its child node, resulting in the encoding [3,2,1] for node $v_2$. Simultaneously, we pass this encoding to its parent node $v_1$, becoming [[3,2,1]];
    \item Following the same steps as mentioned above, we obtain the encoding [4,3,2,1] for the root node $v_1$.
\end{itemize}

\bigskip

\begin{figure}[htbp]
\centering

%\resizebox{\textwidth}{!}{%
\begin{tikzpicture}[
    font=\normalsize,
    vertex/.style={
        circle,
        draw=black,
        line width=0.8pt,
        minimum size=5.5mm,
        inner sep=0pt
    },
    treeedge/.style={
        draw=black,
        line width=0.8pt
    },
    oldedge/.style={
        draw=graytext,
        line width=0.8pt
    },
    steparrow/.style={
        -{Stealth[length=2.5mm,width=1.6mm]},
        draw=arrowblue,
        line width=1pt
    },
    stepnumber/.style={
        circle,
        fill=arrowblue,
        text=white,
        font=\bfseries\scriptsize,
        minimum size=3mm,
        inner sep=0pt
    },
    panel/.style={
        rounded corners=2mm,
        fill=panelblue,
        draw=none
    },
    encoding/.style={
        font=\small,
        text=encodingred,
        anchor=west
    },
    encb/.style={
        font=\small,
        text=arrowblue,
        anchor=west
    },
    oldencoding/.style={
        font=\small,
        text=graytext,
        anchor=west
    },
    vlabel/.style={
        font=\small,
        anchor=east
    }
]

% --------------------------------------------------
% Panel positions and dimensions
% --------------------------------------------------
\def\xA{2}
\def\xB{4.9}
\def\xC{8.4}
\def\xD{11.9}

\def\panelwidth{3.3}
\def\panelbottom{0.35}
\def\paneltop{3.65}

% Vertex coordinates relative to each panel's left edge
\def\xV{0.75}
\def\xRoot{1.35}
\def\xVFour{2.35}

% Vertical positions
\def\yTop{3.0}
\def\yMid{2.0}
\def\yRoot{0.9}

% --------------------------------------------------
% Panel (a): Original tree
% --------------------------------------------------
\fill[panel]
    (\xA+0.3,\panelbottom)
    rectangle
    (\xA+\panelwidth-0.6,\paneltop);

\draw[treeedge]
    (\xA+\xV,\yTop) -- (\xA+\xV,\yMid)
    (\xA+\xV,\yMid) -- (\xA+\xRoot,\yRoot)
    (\xA+\xVFour, \yMid) -- (\xA+\xRoot,\yRoot);

\node[vertex,fill=white]
    at (\xA+\xV,\yTop) {};
\node[vertex,fill=white]
    at (\xA+\xV,\yMid) {};
\node[vertex,fill=rootpurple]
    at (\xA+\xRoot,\yRoot) {};
\node[vertex,fill=white]
    at (\xA+\xVFour,\yMid) {};

\node[vlabel]
    at (\xA+\xV+0.35,\yTop) {$v_1$};
\node[vlabel]
    at (\xA+\xV+0.35,\yMid) {$v_2$};
\node[vlabel]
    at (\xA+\xRoot+0.35,\yRoot) {$v_3$};
\node[vlabel]
    at (\xA+\xVFour+0.35,\yMid) {$v_4$};

% --------------------------------------------------
% Transition 1
% --------------------------------------------------

% --------------------------------------------------
% Panel (b): First encoding
% --------------------------------------------------
\fill[panel]
    (\xB+0.3,\panelbottom)
    rectangle
    (\xB+\panelwidth,\paneltop);

\draw[treeedge]
    (\xB+\xV,\yTop) -- (\xB+\xV,\yMid)
    (\xB+\xV,\yMid) -- (\xB+\xRoot,\yRoot)
    (\xB+\xVFour,\yMid) -- (\xB+\xRoot,\yRoot);

\node[vertex,fill=white]
    at (\xB+\xV,\yTop) {};
\node[vertex,fill=white]
    at (\xB+\xV,\yMid) {};
\node[vertex,fill=rootpurple]
    at (\xB+\xRoot,\yRoot) {};
\node[vertex,fill=white]
    at (\xB+\xVFour,\yMid) {};

\node[vlabel]
    at (\xB+\xV+0.35,\yTop) {$v_1$};
\node[vlabel]
    at (\xB+\xV+0.35,\yMid) {$v_2$};
\node[vlabel]
    at (\xB+\xRoot+0.35,\yRoot) {$v_3$};
\node[vlabel]
    at (\xB+\xVFour+0.35,\yMid) {$v_4$};

\node[encoding]
    at (\xB+0.9,\yTop) {$[1]$};
\node[encb]
    at (\xB+0.9,\yMid) {$[[1]]$};
\node[encoding]
    at (\xB+2.5,\yMid) {$[1]$};
\node[encb]
    at (\xB+1.5,\yRoot) {$[[1]]$};

% --------------------------------------------------
% Transition 2
% --------------------------------------------------

% --------------------------------------------------
% Panel (c): Second encoding
% --------------------------------------------------
\fill[panel]
    (\xC+0.3,\panelbottom)
    rectangle
    (\xC+\panelwidth,\paneltop);

\draw[oldedge]
    (\xC+\xV,\yTop) -- (\xC+\xV,\yMid);
\draw[treeedge]
    (\xC+\xV,\yMid) -- (\xC+\xRoot,\yRoot);
\draw[oldedge]
    (\xC+\xVFour,\yMid) -- (\xC+\xRoot,\yRoot);

\node[vertex,draw=graytext,fill=white]
    at (\xC+\xV,\yTop) {};
\node[vertex,fill=white]
    at (\xC+\xV,\yMid) {};
\node[vertex,fill=rootpurple]
    at (\xC+\xRoot,\yRoot) {};
\node[vertex,draw=graytext,fill=white]
    at (\xC+\xVFour,\yMid) {};

\node[vlabel,text=graytext]
    at (\xC+\xV+0.35,\yTop) {$v_1$};
\node[vlabel]
    at (\xC+\xV+0.35,\yMid) {$v_2$};
\node[vlabel]
    at (\xC+\xRoot+0.35,\yRoot) {$v_3$};
\node[vlabel,text=graytext]
    at (\xC+\xVFour+0.35,\yMid) {$v_4$};

\node[oldencoding]
    at (\xC+0.9,\yTop) {$[1]$};
\node[encoding]
    at (\xC+0.9,\yMid) {$[2,1]$};
\node[oldencoding]
    at (\xC+2.5,\yMid) {$[1]$};
\node[encb]
    at (\xC+1.5,\yRoot) {$[1],[2,1]$};

% --------------------------------------------------
% Transition 3
% --------------------------------------------------

% --------------------------------------------------
% Panel (d): Final encoding
% --------------------------------------------------
\fill[panel]
    (\xD+0.3,\panelbottom)
    rectangle
    (\xD+\panelwidth,\paneltop);

\draw[oldedge]
    (\xD+\xV,\yTop) -- (\xD+\xV,\yMid)
    (\xD+\xV,\yMid) -- (\xD+\xRoot,\yRoot)
    (\xD+\xVFour,\yMid) -- (\xD+\xRoot,\yRoot);

\node[vertex,draw=graytext,fill=white]
    at (\xD+\xV,\yTop) {};
\node[vertex,draw=graytext,fill=white]
    at (\xD+\xV,\yMid) {};
\node[vertex,fill=rootpurple]
    at (\xD+\xRoot,\yRoot) {};
\node[vertex,draw=graytext,fill=white]
    at (\xD+\xVFour,\yMid) {};

\node[vlabel,text=graytext]
    at (\xD+\xV+0.35,\yTop) {$v_1$};
\node[vlabel,text=graytext]
    at (\xD+\xV+0.35,\yMid) {$v_2$};
\node[vlabel]
    at (\xD+\xRoot++0.35,\yRoot) {$v_3$};
\node[vlabel,text=graytext]
    at (\xD+\xVFour+0.35,\yMid) {$v_4$};

\node[oldencoding]
    at (\xD+0.9,\yTop) {$[1]$};
\node[oldencoding]
    at (\xD+0.9,\yMid) {$[2,1]$};
\node[oldencoding]
    at (\xD+2.5,\yMid) {$[1]$};
\node[encoding]
    at (\xD+1.5,\yRoot) {$[4,1,2,1]$};

\draw[steparrow]
    (4.75,1.75) -- (5.15,1.75);
\node[stepnumber]
    at (4.95,2.08) {$1$};

\draw[steparrow]
    (8.25,1.75) -- (8.65,1.75);
\node[stepnumber]
    at (8.45,2.08) {$2$};
    
\draw[steparrow]
    (11.75,1.75) -- (12.15,1.75);
\node[stepnumber]
    at (11.95,2.08) {$3$};

\node[
    font=\scriptsize,
    text=navy,
    align=center,
    text width=2.7cm
] at (\xA+1.4,0.05)
    {(a) Original tree};

\node[
    font=\scriptsize,
    text=navy,
    align=center,
    text width=4cm
] at (\xB+1.8,0.05)
    {(b)Processing $v_1$ and $v_4$};

\node[
    font=\scriptsize,
    text=navy,
    align=center,
    text width=4cm
] at (\xC+1.8,0.05)
    {(c)Processing $v_2$};

\node[
    font=\scriptsize,
    text=navy,
    align=center,
    text width=2.7cm
] at (\xD+1.8,0.05)
    {(d) Final encoding};
\end{tikzpicture}%
%}

\caption{Step-by-step encoding of another rooted tree with four nodes.}
\label{fig-02}
\end{figure}
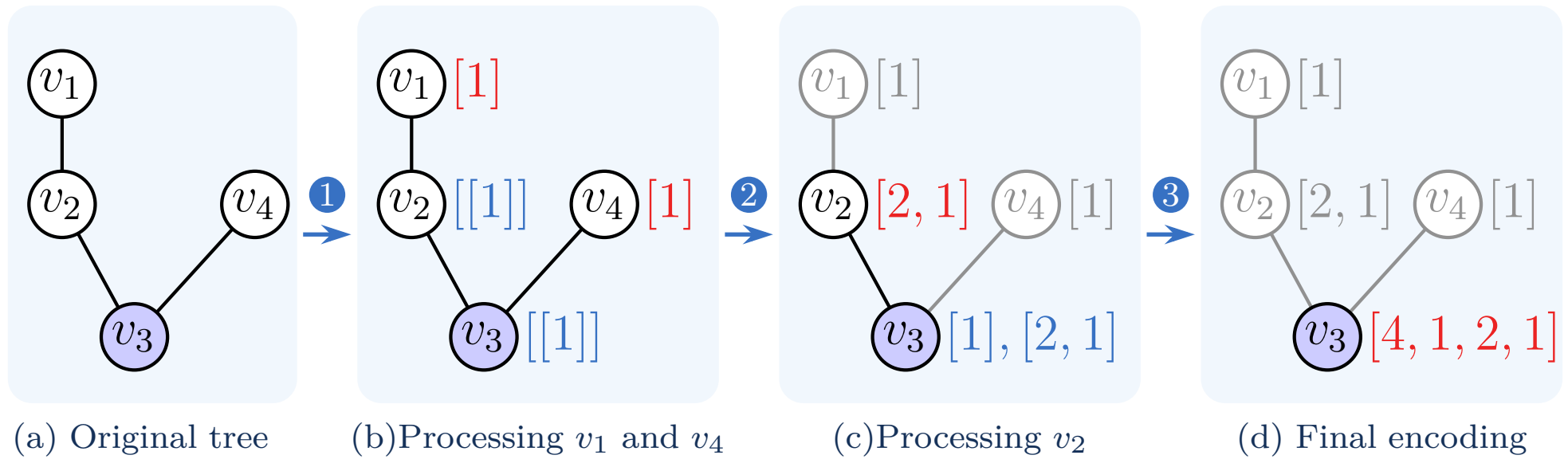

\begin{example}
As shown in Figure \ref{fig-02}, we have another four-node tree where $v_3$ is the root node. The encoding process starts from the leaf nodes $v_4$ and $v_1$. The encoding process is as follows,
\end{example}

\begin{itemize}
    \item The encoding of both $v_1$ and $v_4$ is [1], and this encoding is passed to their respective parent nodes $v_2$ and $v_3$, resulting in [[1]] for both;
    \item Considering the removal of the previously encoded subtree rooted at nodes $v_1$ and $v_4$, we now select a new leaf node $v_2$. The total number of nodes in the subtree rooted at $v_2$ is 2. We then concatenate the encoding [1] of its child node $v_1$, resulting in the encoding [2,1] for node $v_2$. It should be noted that although $v_3$ is also a new leaf node, it has been designated as the root node, so no final encoding is performed at this point. Instead, the encoding [2,1] received from $v_2$ as a child node is temporarily stored together with the previously stored temporary encoding [1], resulting in a new temporary encoding for $v_3$ as [[1],[2,1]].
    \item Following the same steps described above, we continue by removing the leaf nodes from the previous step, leaving only the root node $v_3$. Its encoding is then obtained as [4,1,2,1], where the leading 4 is the total number of nodes in the subtree rooted at  $v_3$, and the remaining 1, 2, and 1 are the encodings of its child nodes, [1] and [2,1], sorted and concatenated after stripping the square brackets.
\end{itemize}

    From Figures \ref{fig-01} and \ref{fig-02}, it can be observed that although they are structurally identical four-node trees, their encodings are different due to the specification of different root nodes. The reason for this lies in the fact that along different root nodes, the subtrees of each node will vary, resulting in different encodings.

\subsection{Encoding of Unrooted Trees}
For an unrooted tree $T(V, E)$ with a unique center node $c$, we transform it into a rooted tree $T'(V, E, c)$. If $T$ has a bicentral structure (a two-node center $\{c_1, c_2\}$), we remove the edge connecting these two nodes to decompose $T$ into two rooted subtrees $T'_i(V_i, E_i, c_i)$ for $i \in \{1, 2\}$, where the roots $c_i$ correspond to the original center nodes.

\begin{defn}
\label{defn210}
     Let $T=(V,E)$ be an unrooted tree with $n$ nodes. Its encoding is defined as follows:

{\footnotesize
    \begin{equation}
    \label{08_1}
    \text{code}(T) =
        \begin{cases} 
            \text{code}(T'(V,E,c)) & \text{The center of } T \text{ is a single-node } c, \\[15pt] % 增加行间距
             [op(\text{code}(T'_1)), op(\text{code}(T'_2))] & 
                \parbox[t]{0.45\textwidth}{ % 限制宽度以解决 Overfull \hbox
                The center of $T$ consists of two nodes, $c_1$ and $c_2$. 
                The two subtrees $T_1$ and $T_2$, with $c_1$ and $c_2$ 
                as their root nodes, are encoded in ascending order 
                based on dictionary sequence.
         }
        \end{cases}
    \end{equation}
}
Here $T'(V,E,c)$, $T'_1(V_1,E_1,c_1)$, and $T'_2(V_2,E_2,c_2)$ are derived from the original tree $T$, as discussed in the preceding definitions. $op(code(T'_i(V_i,E_i,c_i)))$ $(i=1,2)$ represents extracting the sequence of integers from the encoding of $T'_i(V_i,E_i,c_i)$ $($i.e., removing a pair of square brackets $[\; ]$ from the encoding of $T'_i(V_i,E_i,c_i))$.
\end{defn}

\begin{figure}[H]
\centering
\begin{tikzpicture}[
    font=\normalsize,
    vertex/.style={
        circle,
        draw=black,
        line width=0.8pt,
        minimum size=5.5mm,
        inner sep=0pt,
        fill=white
    },
    vertexgray/.style={
        circle,
        draw=graytext,
        line width=0.8pt,
        minimum size=5.5mm,
        inner sep=0pt,
        fill=white
    },
    treeedge/.style={
        draw=black,
        line width=0.8pt
    },
    oldedge/.style={
        draw=graytext,
        line width=0.8pt
    },
    steparrow/.style={
        -{Stealth[length=2.5mm,width=1.6mm]},
        draw=arrowblue,
        line width=1pt
    },
    stepnumber/.style={
        circle,
        fill=arrowblue,
        text=white,
        font=\bfseries\scriptsize,
        minimum size=3.2mm,
        inner sep=0pt
    },
    panel/.style={
        rounded corners=2mm,
        fill=panelblue,
        draw=none
    },
    encoding/.style={
        font=\small,
        text=encodingred,
        anchor=west
    },
    oldencoding/.style={
        font=\small,
        text=graytext,
        anchor=west
    },
    encb/.style={font=\small,
         text=arrowblue,
        anchor=west
    },
    vlabel/.style={
        font=\small,
        anchor=east
    }
]

% --------------------------------------------------
% Color definitions (if not already defined in preamble)
% --------------------------------------------------
\definecolor{arrowblue}{RGB}{30,100,200}
\definecolor{panelblue}{RGB}{235,242,250}
\definecolor{encodingred}{RGB}{200,30,30}
\definecolor{graytext}{RGB}{160,160,160}

% --------------------------------------------------
% Panel positions and dimensions
% --------------------------------------------------
\def\xA{1.0}
\def\xB{4.0}
\def\xC{7.75}
\def\xD{11.3}
\def\panelwidth{2.9}
\def\panelbottom{-0.35}
\def\paneltop{4.05}

% Vertex x offset inside each panel
\def\xV{1.0}

% Vertical positions (bottom to top)
\def\yOne{0.5}
\def\yTwo{1.5}
\def\yThree{2.5}
\def\yFour{3.5}

% --------------------------------------------------
% Panel (a): Original tree
% --------------------------------------------------
\fill[panel]
    (\xA,\panelbottom) rectangle (\xA+\panelwidth-0.5,\paneltop);

\draw[treeedge]
    (\xA+\xV+0.2,\yFour) -- (\xA+\xV+0.2,\yThree)
    (\xA+\xV+0.2,\yThree) -- (\xA+\xV+0.2,\yTwo)
    (\xA+\xV+0.2,\yTwo) -- (\xA+\xV+0.2,\yOne);

\node[vertex] at (\xA+\xV+0.2,\yFour) {};
\node[vertex] at (\xA+\xV+0.2,\yThree) {};
\node[vertex] at (\xA+\xV+0.2,\yTwo) {};
\node[vertex] at (\xA+\xV+0.2,\yOne) {};

\node[vlabel] at (\xA+\xV+0.2-0.35,\yFour) {$v_4$};
\node[vlabel] at (\xA+\xV+0.2-0.35,\yThree) {$v_3$};
\node[vlabel] at (\xA+\xV+0.2-0.35,\yTwo) {$v_2$};
\node[vlabel] at (\xA+\xV+0.2-0.35,\yOne) {$v_1$};

% --------------------------------------------------
% Transition ①
% --------------------------------------------------
\draw[steparrow] (\xA+\panelwidth-0.5,2.0) -- (\xB-0.2,2.0);
\node[stepnumber] at (\xA+\panelwidth-0.3,2.35) {1};

% --------------------------------------------------
% Panel (b): After first encoding
% --------------------------------------------------
\fill[panel]
    (\xB-0.2,\panelbottom) rectangle (\xB+\panelwidth+0.25,\paneltop);

\draw[treeedge]
    (\xB+\xV,\yFour) -- (\xB+\xV,\yThree)
    (\xB+\xV,\yThree) -- (\xB+\xV,\yTwo)
    (\xB+\xV,\yTwo) -- (\xB+\xV,\yOne);

\node[vertex] at (\xB+\xV,\yFour) {};
\node[vertex] at (\xB+\xV,\yThree) {};
\node[vertex] at (\xB+\xV,\yTwo) {};
\node[vertex] at (\xB+\xV,\yOne) {};

\node[vlabel] at (\xB+\xV-0.35,\yFour) {$v_4$};
\node[vlabel] at (\xB+\xV-0.35,\yThree) {$v_3$};
\node[vlabel] at (\xB+\xV-0.35,\yTwo) {$v_2$};
\node[vlabel] at (\xB+\xV-0.35,\yOne) {$v_1$};

\node[encoding] at (\xB+\xV+0.35,\yFour) {$[1]$};
\node[encb] at (\xB+\xV+0.35,\yThree) {$[[1]]$};
\node[encb] at (\xB+\xV+0.35,\yTwo) {$[[1]]$};
\node[encoding] at (\xB+\xV+0.35,\yOne) {$[1]$};

% --------------------------------------------------
% Transition ②
% --------------------------------------------------
\draw[steparrow] (\xB+\panelwidth+0.25,2.0) -- (\xC-0.2,2.0);
\node[stepnumber] at (\xB+\panelwidth+0.45,2.35) {2};

% --------------------------------------------------
% Panel (c): After processing edge e_{2,3}
% --------------------------------------------------
\fill[panel]
    (\xC-0.2,\panelbottom) rectangle (\xC+\panelwidth+0.25,\paneltop);

\draw[oldedge]
    (\xC+\xV,\yFour) -- (\xC+\xV,\yThree);
\draw[treeedge]
    (\xC+\xV,\yThree) -- (\xC+\xV,\yTwo);
\draw[oldedge]
    (\xC+\xV,\yTwo) -- (\xC+\xV,\yOne);

\node[vertexgray] at (\xC+\xV,\yFour) {};
\node[vertex] at (\xC+\xV,\yThree) {};
\node[vertex] at (\xC+\xV,\yTwo) {};
\node[vertexgray] at (\xC+\xV,\yOne) {};

\node[vlabel,text=graytext] at (\xC+\xV-0.35,\yFour) {$v_4$};
\node[vlabel] at (\xC+\xV-0.35,\yThree) {$v_3$};
\node[vlabel] at (\xC+\xV-0.35,\yTwo) {$v_2$};
\node[vlabel,text=graytext] at (\xC+\xV-0.35,\yOne) {$v_1$};

\node[oldencoding] at (\xC+\xV+0.35,\yFour) {$[1]$};
\node[encoding] at (\xC+\xV+0.35,\yThree) {$[2,1]$};
\node[encoding] at (\xC+\xV+0.35,\yTwo) {$[2,1]$};
\node[oldencoding] at (\xC+\xV+0.35,\yOne) {$[1]$};

% edge label
\node[font=\scriptsize,text=black,anchor=west] at (\xC+\xV,2.0) {$e_{2,3}$};

% --------------------------------------------------
% Transition ③
% --------------------------------------------------
\draw[steparrow] (\xC+\panelwidth+0.25,2.0) -- (\xD,2.0);
\node[stepnumber] at (\xC+\panelwidth+0.45,2.35) {3};

% --------------------------------------------------
% Panel (d): Final encoding
% --------------------------------------------------
\fill[panel]
    (\xD,\panelbottom) rectangle (\xD+\panelwidth+2.6,\paneltop);

\draw[oldedge]
    (\xD+\xV,\yFour) -- (\xD+\xV,\yThree);
%\draw[treeedge]
%    (\xD+\xV,\yThree) -- (\xD+\xV,\yTwo);
\draw[oldedge]
    (\xD+\xV,\yTwo) -- (\xD+\xV,\yOne);

\node[vertexgray] at (\xD+\xV,\yFour) {};
\node[vertexgray] at (\xD+\xV,\yThree) {};
\node[vertexgray] at (\xD+\xV,\yTwo) {};
\node[vertexgray] at (\xD+\xV,\yOne) {};

\node[vlabel,text=graytext] at (\xD+\xV-0.35,\yFour) {$v_4$};
\node[vlabel,text=graytext] at (\xD+\xV-0.35,\yThree) {$v_3$};
\node[vlabel,text=graytext] at (\xD+\xV-0.35,\yTwo) {$v_2$};
\node[vlabel,text=graytext] at (\xD+\xV-0.35,\yOne) {$v_1$};

\node[oldencoding] at (\xD+\xV+0.35,\yFour) {$[1]$};
\node[oldencoding] at (\xD+\xV+0.35,\yThree) {$[2,1]$};
\node[oldencoding] at (\xD+\xV+0.35,\yTwo) {$[2,1]$};
\node[oldencoding] at (\xD+\xV+0.35,\yOne) {$[1]$};

% final code
\node[font=\footnotesize,text=encodingred,anchor=west]
    at (\xD+\xV+1.15,2.0) {$code(T)=[2,1,2,1]$};

\node[
    font=\scriptsize,
    text=navy,
    align=center,
    text width=2.7cm
] at (\xA+1.15,-0.1)
    {(a) Original tree};

\node[
    font=\scriptsize,
    text=navy,
    align=center,
    text width=4cm
] at (\xB+1.5,-0.1)
    {(b)Processing $v_4$ and $v_1$};

\node[
    font=\scriptsize,
    text=navy,
    align=center,
    text width=4cm
] at (\xC+1.5,-0.1)
    {(c)Processing $v_3$ and $v_2$};

\node[
    font=\scriptsize,
    text=navy,
    align=center,
    text width=2.7cm
] at (\xD+1.8,-0.1)
    {(d) Final encoding};
\end{tikzpicture}
\caption{Trees with four nodes}
\label{Fig-03}
\end{figure}
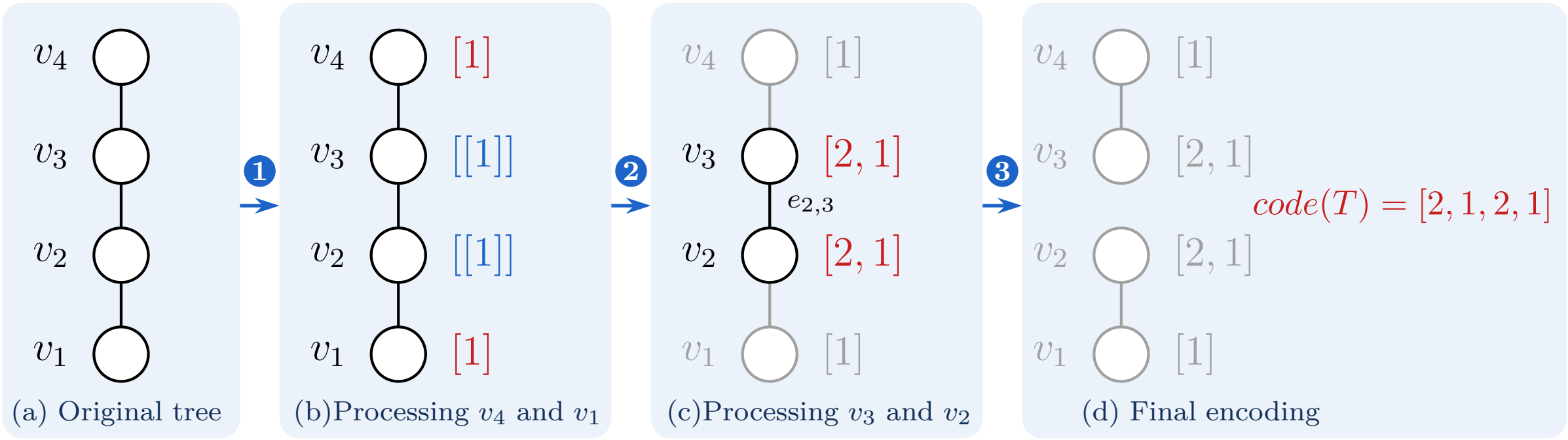

\begin{example}
Figure \ref{Fig-03} illustrates an unrooted four-node tree, where the encoding process initiates from leaf nodes $v_4$ and $v_1$ through the following steps:
\end{example}
\begin{itemize}
    \item Leaf-to-Parent Propagation: The initial encodings for leaf nodes $v_1$ and $v_4$ are both $[1]$. These encodings are propagated to their respective parent nodes, $v_2$ and $v_3$, resulting in the intermediate string $[[1]]$ for each.
    \item Identification of the Graph Center: After pruning the encoded leaves $v_1$ and $v_4$, the remaining structure contains exactly two nodes, namely $v_2$ and $v_3$. Hence the center of the original tree is identified as the node pair $(v_2, v_3)$, which now act as the leaves of the surviving subtree. Following the procedure in the previous example, both $v_2$ and $v_3$ yield the encoding $[2,1]$. Note that the edge $e_{2,3}$ between them is excluded here, as required by the definition of unrooted tree encoding.
    \item Canonical Concatenation: The encodings of the two central nodes are sorted in ascending (lexicographical) order, yielding $[[2,1], [2,1]]$. By removing the outer brackets of each individual encoding and concatenating the sequences, we derive the string $2,1,2,1$. Thus, the final canonical encoding for the original unrooted tree is $\text{code}(T) = [2,1,2,1]$.
\end{itemize}

For single-center unrooted trees, the process is relatively straightforward and precludes the need for additional examples. The method employs an iterative pruning strategy: leaf nodes of the subtree are successively encoded and removed until only the central structure remains. This procedure, reminiscent of the classical algorithm for locating a tree's center, naturally converges to one or two remaining nodes. By integrating center identification directly into the encoding sequence, the requirement for pre-identifying centers is eliminated, thereby optimizing the algorithm's computational efficiency.

\subsection{Encoding Unrooted Colored Trees}
We now extend the encoding approach to unrooted trees whose nodes are colored black or white. Unlike the uncolored trees considered above, these trees must be distinguished by both their structures and their node colors. To develop an encoding for unrooted colored trees, we first define the encoding of rooted colored trees and then extend it to the unrooted case.
\begin{defn}
Let $T=(V,E,r)$ be a rooted tree whose nodes are colored black or white,  where $V$ is the set of nodes, $E$ is the set of edges, and $r$ is the root node of the tree. For any node $v \in V$, its encoding is defined as the following sequence of integers,

\begin{eqnarray}
     \label{08_2}
    code(v)=\left\{
    \begin{array}{ll}
        [sign(v) \times size(v),op(code(w_1)),\cdots,op(code(w _k))] &  v \mbox{ is non-leaf node}, \\
        \left[ sign(v) \times 1 \right] & v\mbox{ is leaf node},
    \end{array}
    \right.
\end{eqnarray}
where
\begin{eqnarray*}
     \label{09}
    sign(v)=\left\{
    \begin{array}{ll}
        -1  &  v \mbox{ is black node}, \\
        1 &  v\mbox{ is white node},
    \end{array}
    \right.
\end{eqnarray*}
Where $size(v)$, $w_1,\cdots,w_k$, and $op(code(w_j)) (j=1, \cdots, k)$ are defined in \textnormal{Definition \ref{def01}}.
\end{defn}

\begin{figure}[H]
\centering
\begin{tikzpicture}[
    font=\normalsize,
    vertex/.style={
        circle, draw=black, line width=0.75pt,
        minimum size=5.4mm, inner sep=0pt, fill=white
    },
    vertexblack/.style={
        circle, draw=black, line width=0.75pt,
        minimum size=5.4mm, inner sep=0pt, fill=black
    },
    vertexgray/.style={
        circle, draw=graytext, line width=0.75pt,
        minimum size=5.4mm, inner sep=0pt, fill=gray!35
    },
    vertexgrayw/.style={
        circle, draw=graytext, line width=0.75pt,
        minimum size=5.4mm, inner sep=0pt, fill=white
    },
    treeedge/.style={draw=black, line width=0.75pt},
    oldedge/.style={draw=graytext, line width=0.75pt},    
    steparrow/.style={
        -{Stealth[length=2mm,width=1.5mm]},
        draw=arrowblue, line width=1.05pt
    },
    stepnumber/.style={
        circle, fill=arrowblue, text=white,
        font=\bfseries\footnotesize, minimum size=4.0mm, inner sep=0pt
    },
    panel/.style={
        rounded corners=2.2mm, fill=panelblue, draw=none
    },
    enc/.style={font=\scriptsize, text=encodingred},
    encb/.style={font=\scriptsize, text=arrowblue},
    oldenc/.style={font=\scriptsize, text=graytext},
    vl/.style={font=\footnotesize}
]

\definecolor{arrowblue}{RGB}{30,100,200}
\definecolor{panelblue}{RGB}{235,242,250}
\definecolor{encodingred}{RGB}{200,30,30}
\definecolor{graytext}{RGB}{150,150,150}

% ============================================================
% 布局参数（可根据需要微调）
% ============================================================
\def\pw{6.4}          % panel 宽度
\def\gap{1.3}         % 同行两个 panel 间距

% 第1行
\def\yA{9.8}
\def\xA{0.0}
\pgfmathsetmacro{\xB}{\xA+\pw+\gap}

% 第2行
\def\yC{3.0}
\def\xC{0.0}
\pgfmathsetmacro{\xD}{\xC+\pw+\gap}

% 第3行（居中）
\def\yE{-3.9}
\def\xE{0.0}
% \pgfmathsetmacro{\xE}{(\xA+\xB+\pw)/2 - \pw/2}

% 树内部相对坐标（以每个 panel 左下角附近为参考）
% v6 top → v9 bottom
\newcommand{\vx}{2.9}          % 树中心 x 偏移
\newcommand{\vy}{0.35}         % 树底部 y 偏移

% ============================================================
% Panel A：原始树
% ============================================================
\fill[panel] (\xA-0.55,\yA-0.45) rectangle (\xA+\pw+0.2,\yA+6);

\begin{scope}[shift={(\xA+\vx,\yA+\vy)}]
  \draw[treeedge]
    (0,4.9) -- (0,3.9)           % v6-v5
    (0,3.9) -- (-1.65,2.55)      % v5-v3
    (0,3.9) -- (1.65,2.55)       % v5-v8
    (-1.65,2.55) -- (-2.55,1.25) % v3-v1
    (-1.65,2.55) -- (-0.45,1.25) % v3-v4
    (-0.45,1.25) -- (-1.35,0.05) % v4-v2
    (-0.45,1.25) -- (0.45,0.05)  % v4-v7
    (1.65,2.55) -- (2.75,2.55)   % v8-v11
    (1.65,2.55) -- (1.65,1.25)   % v8-v10
    (1.65,1.25) -- (1.65,0.05);  % v10-v9

  \node[vertexblack] at (0,4.9) {};      % v6
  \node[vertex]      at (0,3.9) {};      % v5
  \node[vertex]      at (-1.65,2.55) {}; % v3
  \node[vertex]      at (1.65,2.55) {};  % v8
  \node[vertex]      at (-2.55,1.25) {}; % v1
  \node[vertexblack] at (-0.45,1.25) {}; % v4
  \node[vertexblack] at (-1.35,0.05) {}; % v2
  \node[vertex]      at (0.45,0.05) {};  % v7
  \node[vertex]      at (2.75,2.55) {};  % v11
  \node[vertex]      at (1.65,1.25) {};  % v10
  \node[vertex]      at (1.65,0.05) {};  % v9

  \node[vl,left]  at (-0.2,4.9) {$v_6$};
  \node[vl,left]  at (-0.2,3.9) {$v_5$};
  \node[vl,left]  at (-1.85,2.55) {$v_3$};
  \node[vl,left]  at (1.45,2.55) {$v_8$};
  \node[vl,left]  at (-2.75,1.25) {$v_1$};
  \node[vl,left]  at (-0.65,1.25) {$v_4$};
  \node[vl,left]  at (-1.55,0.05) {$v_2$};
  \node[vl,left]  at (0.25,0.05) {$v_7$};
  \node[vl,right] at (2.95,2.55) {$v_{11}$};
  \node[vl,left]  at (1.45,1.25) {$v_{10}$};
  \node[vl,left]  at (1.45,0.05) {$v_9$};
  
  \node[
    font=\scriptsize,
    text=navy,
    align=center,
    text width=2.7cm
] at (0,-1)
    {(a) Original tree};

\end{scope}

% Arrow ①
\draw[steparrow] (\xA+\pw+0.2,\yA+2.7) -- (\xB-0.6,\yA+2.7);
\node[stepnumber] at (\xA+\pw+0.45,\yA+3.15) {1};

% ============================================================
% Panel B：步骤①后
% ============================================================
\fill[panel] (\xB-0.55,\yA-0.45) rectangle (\xB+\pw+0.35,\yA+6);

\begin{scope}[shift={(\xB+\vx,\yA+\vy)}]
  \draw[treeedge]
    (0,4.9) -- (0,3.9)
    (0,3.9) -- (-1.65,2.55)
    (0,3.9) -- (1.65,2.55)
    (-1.65,2.55) -- (-2.55,1.25)
    (-1.65,2.55) -- (-0.45,1.25)
    (-0.45,1.25) -- (-1.35,0.05)
    (-0.45,1.25) -- (0.45,0.05)
    (1.65,2.55) -- (2.75,2.55)
    (1.65,2.55) -- (1.65,1.25)
    (1.65,1.25) -- (1.65,0.05);

  \node[vertexblack] at (0,4.9) {};
  \node[vertex]      at (0,3.9) {};
  \node[vertex]      at (-1.65,2.55) {};
  \node[vertex]      at (1.65,2.55) {};
  \node[vertex]      at (-2.55,1.25) {};
  \node[vertexblack] at (-0.45,1.25) {};
  \node[vertexblack] at (-1.35,0.05) {};
  \node[vertex]      at (0.45,0.05) {};
  \node[vertex]      at (2.75,2.55) {};
  \node[vertex]      at (1.65,1.25) {};
  \node[vertex]      at (1.65,0.05) {};

  \node[vl,left]  at (-0.2,4.9) {$v_6$};
  \node[vl,left]  at (-0.2,3.9) {$v_5$};
  \node[vl,left]  at (-1.85,2.55) {$v_3$};
  \node[vl,left]  at (1.45,2.55) {$v_8$};
  \node[vl,left]  at (-2.75,1.25) {$v_1$};
  \node[vl,left]  at (-0.65,1.25) {$v_4$};
  \node[vl,left]  at (-1.55,0.05) {$v_2$};
  \node[vl,left]  at (0.25,0.05) {$v_7$};
  \node[vl,right] at (2.95,2.55) {$v_{11}$};
  \node[vl,left]  at (1.45,1.25) {$v_{10}$};
  \node[vl,left]  at (1.45,0.05) {$v_9$};

  % encodings
  \node[enc] at (0.5,5.3) {$[-1]$};
  \node[encb] at (0.6,4.4) {$[[-1]]$};
  \node[encb] at (-1.7,3.1) {$[[1]]$};
  \node[encb] at (1.6,3.1) {$[[1]]$};
  \node[enc] at (-2.9,0.8) {$[1]$};
  \node[encb] at (-0.2,1.8) {$[[-1],[1]]$};
  \node[enc] at (-1.7,-0.4) {$[-1]$};
  \node[enc] at (0.8,-0.4) {$[1]$};
  \node[enc] at (3.1,3.1) {$[1]$};
  \node[encb] at (2.0,1.8) {$[[1]]$};
  \node[enc] at (2.0,-0.4) {$[1]$};
  
    \node[
    font=\scriptsize,
    text=navy,
    align=center,
    text width=6cm
] at (0,-1)
    {(b) Processing leaf node$v_6,v_1,v_2,v_7,v_9$ and $v_{11}$};
\end{scope}

% ============================================================
% Arrow ②：与步骤①、③同为水平短箭头，置于第二幅子图右侧
% ============================================================
\draw[steparrow] (\xB+6.75,\yA+2.7) -- (\xB+7.25,\yA+2.7);
\node[stepnumber] at (\xB+7.0,\yA+3.15) {2};

% ============================================================
% Panel C：步骤②后
% ============================================================
\fill[panel] (\xC-0.55,\yC-0.45) rectangle (\xC+\pw+0.2,\yC+6);

\begin{scope}[shift={(\xC+\vx,\yC+\vy)}]
  \draw[oldedge]
    (0,4.9) -- (0,3.9)
    (-1.65,2.55) -- (-2.55,1.25)
    (-0.45,1.25) -- (-1.35,0.05)
    (-0.45,1.25) -- (0.45,0.05)    
    (1.65,1.25) -- (1.65,0.05)
    (1.65,2.55) -- (2.75,2.55);
  \draw[treeedge]
    (0,3.9) -- (-1.65,2.55)
    (0,3.9) -- (1.65,2.55)
    (-1.65,2.55) -- (-0.45,1.25)  
    (1.65,2.55) -- (1.65,1.25);

  \node[vertexgray]  at (0,4.9) {};
  \node[vertex]      at (0,3.9) {};
  \node[vertex]      at (-1.65,2.55) {};
  \node[vertex]      at (1.65,2.55) {};
  \node[vertexgrayw] at (-2.55,1.25) {};
  \node[vertexblack] at (-0.45,1.25) {};
  \node[vertexgray]  at (-1.35,0.05) {};
  \node[vertexgrayw] at (0.45,0.05) {};
  \node[vertexgrayw] at (2.75,2.55) {};
  \node[vertex]      at (1.65,1.25) {};
  \node[vertexgrayw] at (1.65,0.05) {};

  \node[vl,left,text=graytext]  at (-0.2,4.9) {$v_6$};
  \node[vl,left]  at (-0.2,3.9) {$v_5$};
  \node[vl,left]  at (-1.85,2.55) {$v_3$};
  \node[vl,left]  at (1.45,2.55) {$v_8$};
  \node[vl,left,text=graytext]  at (-2.75,1.25) {$v_1$};
  \node[vl,left]  at (-0.65,1.25) {$v_4$};
  \node[vl,left,text=graytext]  at (-1.55,0.05) {$v_2$};
  \node[vl,left,text=graytext]  at (0.25,0.05) {$v_7$};
  \node[vl,right,text=graytext] at (2.95,2.55) {$v_{11}$};
  \node[vl,left]  at (1.45,1.25) {$v_{10}$};
  \node[vl,left,text=graytext]  at (1.45,0.05) {$v_9$};

  \node[oldenc] at (0.5,5.3) {$[-1]$};
  \node[encb] at (0.6,4.4) {$[[-1]]$};
  % \node[encb] at (-1.7,3.1) {$[[1]]$};
  % \node[encb] at (1.6,3.1) {$[[1]]$};
  \node[oldenc] at (-2.9,0.8) {$[1]$};
  % \node[encb] at (-0.2,1.8) {$[[-1],[1]]$};
  \node[oldenc] at (-1.7,-0.4) {$[-1]$};
  \node[oldenc] at (0.8,-0.4) {$[1]$};
  \node[oldenc] at (3.1,3.1) {$[1]$};
  %\node[encb] at (2.0,1.8) {$[[1]]$};
  \node[oldenc] at (2.0,-0.4) {$[1]$};

  % \node[oldenc] at (0.4,5.1) {$[-1]$};
  \node[encb] at (0.6,4.4) {$[[-1]]$};
  \node[encb] at (-2.3,3.1) {$[[1],[-3,-1,1]]$};
  \node[encb] at (1.9,3.1) {$[[1],[2,1]]$};
  %\node[oldenc] at (-3.1,1.0) {$[1]$};
  \node[enc] at (-0.1,1.8) {$[-3,-1,1]$};
  %\node[oldenc] at (-1.7,-0.25) {$[-1]$};
  %\node[oldenc] at (0.7,-0.25) {$[1]$};
  %\node[oldenc] at (3.1,2.4) {$[1]$};
  \node[enc] at (2.0,1.8) {$[2,1]$};
  %\node[oldenc] at (2.0,-0.25) {$[1]$};
    
   \node[
    font=\scriptsize,
    text=navy,
    align=center,
    text width=6cm
] at (0,-1)
    {(c) Processing $v_4$ and $v_{10}$};
\end{scope}

% Arrow ③
\draw[steparrow] (\xC+\pw+0.2,\yC+2.7) -- (\xD-0.6,\yC+2.7);
\node[stepnumber] at (\xC+\pw+0.45,\yC+3.15) {3};

% ============================================================
% Panel D：步骤③后
% ============================================================
\fill[panel] (\xD-0.55,\yC-0.45) rectangle (\xD+\pw+0.35,\yC+6);

\begin{scope}[shift={(\xD+\vx,\yC+\vy)}]
  \draw[oldedge]
    (0,4.9) -- (0,3.9)
    (-1.65,2.55) -- (-2.55,1.25)
    (-0.45,1.25) -- (-1.35,0.05)
    (-0.45,1.25) -- (0.45,0.05)
    (1.65,2.55) -- (2.75,2.55)
    (1.65,1.25) -- (1.65,0.05)    
    (-1.65,2.55) -- (-0.45,1.25)
    (1.65,2.55) -- (1.65,1.25);
  \draw[treeedge]
    (0,3.9) -- (-1.65,2.55)
    (0,3.9) -- (1.65,2.55);

  \node[vertexgray]  at (0,4.9) {};
  \node[vertex]      at (0,3.9) {};
  \node[vertex]      at (-1.65,2.55) {};
  \node[vertex]      at (1.65,2.55) {};
  \node[vertexgrayw] at (-2.55,1.25) {};
  \node[vertexgray]  at (-0.45,1.25) {};
  \node[vertexgray]  at (-1.35,0.05) {};
  \node[vertexgrayw] at (0.45,0.05) {};
  \node[vertexgrayw] at (2.75,2.55) {};
  \node[vertexgrayw] at (1.65,1.25) {};
  \node[vertexgrayw] at (1.65,0.05) {};

  \node[vl,left,text=graytext]  at (-0.2,4.9) {$v_6$};
  \node[vl,left]  at (-0.2,3.9) {$v_5$};
  \node[vl,left]  at (-1.85,2.55) {$v_3$};
  \node[vl,left]  at (1.45,2.55) {$v_8$};
  \node[vl,left,text=graytext]  at (-2.75,1.25) {$v_1$};
  \node[vl,left,text=graytext]  at (-0.65,1.25) {$v_4$};
  \node[vl,left,text=graytext]  at (-1.55,0.05) {$v_2$};
  \node[vl,left,text=graytext]  at (0.25,0.05) {$v_7$};
  \node[vl,right,text=graytext] at (2.95,2.55) {$v_{11}$};
  \node[vl,left,text=graytext]  at (1.45,1.25) {$v_{10}$};
  \node[vl,left,text=graytext]  at (1.45,0.05) {$v_9$};

  \node[oldenc] at (0.5,5.3) {$[-1]$};
  %\node[encb] at (0.6,4.4) {$[[-1]]$};
  % \node[encb] at (-1.7,3.1) {$[[1]]$};
  % \node[encb] at (1.6,3.1) {$[[1]]$};
  \node[oldenc] at (-2.9,0.8) {$[1]$};
  % \node[encb] at (-0.2,1.8) {$[[-1],[1]]$};
  \node[oldenc] at (-1.7,-0.4) {$[-1]$};
  \node[oldenc] at (0.8,-0.4) {$[1]$};
  \node[oldenc] at (3.1,3.1) {$[1]$};
  %\node[encb] at (2.0,1.8) {$[[1]]$};
  \node[oldenc] at (2.0,-0.4) {$[1]$};

  % \node[oldenc] at (0.4,5.1) {$[-1]$};
  % \node[encb] at (0.6,4.4) {$[[-1]]$};
  %\node[encb] at (-2.3,3.1) {$[[1],[-3,-1,1]]$};
  %\node[encb] at (1.9,3.1) {$[[1],[2,1]]$};
  %\node[oldenc] at (-3.1,1.0) {$[1]$};
  \node[oldenc] at (-0.1,1.8) {$[-3,-1,1]$};
  %\node[oldenc] at (-1.7,-0.25) {$[-1]$};
  %\node[oldenc] at (0.7,-0.25) {$[1]$};
  %\node[oldenc] at (3.1,2.4) {$[1]$};
  \node[oldenc] at (2.0,1.8) {$[2,1]$};
  %\node[oldenc] at (2.0,-0.25) {$[1]$};

  %\node[oldenc] at (0.4,5.1) {$[-1]$};
  \node[encb] at (0,4.4) {$[[-1],[5,-3,-1,1,1],[4,1,2,1]]$};
  \node[enc] at (-2.3,3.1) {$[5,-3,-1,1,1]$};
  \node[enc] at (1.9,3.1) {$[4,1,2,1]$};
  %\node[oldenc] at (-3.1,1.0) {$[1]$};
  %\node[oldenc] at (0.1,1.0) {$[-3,-1,1]$};
  %\node[oldenc] at (-1.7,-0.25) {$[-1]$};
  %\node[oldenc] at (0.7,-0.25) {$[1]$};
  %\node[oldenc] at (3.1,2.4) {$[1]$};
  %\node[oldenc] at (2.05,1.0) {$[2,1]$};
  %\node[oldenc] at (2.0,-0.25) {$[1]$};

  \node[
    font=\scriptsize,
    text=navy,
    align=center,
    text width=6cm
] at (0,-1)
    {(d) Processing $v_3$ and $v_8$};
\end{scope}

% ============================================================
% Arrow ④：与步骤①、③同为水平短箭头，置于第四幅子图右侧
% ============================================================
\draw[steparrow] (\xD+6.75,\yC+2.7) -- (\xD+7.25,\yC+2.7);
\node[stepnumber] at (\xD+7.0,\yC+3.15) {4};

% ============================================================
% Panel E：最终结果（居中）
% ============================================================
\fill[panel] (\xE-0.55,\yE-0.45) rectangle (\xE+\pw+0.35,\yE+6);

\begin{scope}[shift={(\xE+\vx,\yE+\vy)}]
  \draw[oldedge]
    (0,4.9) -- (0,3.9)
    (0,3.9) -- (-1.65,2.55)
    (0,3.9) -- (1.65,2.55)
    (-1.65,2.55) -- (-2.55,1.25)
    (-1.65,2.55) -- (-0.45,1.25)
    (-0.45,1.25) -- (-1.35,0.05)
    (-0.45,1.25) -- (0.45,0.05)
    (1.65,2.55) -- (2.75,2.55)
    (1.65,2.55) -- (1.65,1.25)
    (1.65,1.25) -- (1.65,0.05);

  \node[vertexgray]  at (0,4.9) {};
  \node[vertex]      at (0,3.9) {};
  \node[vertexgrayw] at (-1.65,2.55) {};
  \node[vertexgrayw] at (1.65,2.55) {};
  \node[vertexgrayw] at (-2.55,1.25) {};
  \node[vertexgray]  at (-0.45,1.25) {};
  \node[vertexgray]  at (-1.35,0.05) {};
  \node[vertexgrayw] at (0.45,0.05) {};
  \node[vertexgrayw] at (2.75,2.55) {};
  \node[vertexgrayw] at (1.65,1.25) {};
  \node[vertexgrayw] at (1.65,0.05) {};

  \node[vl,left,text=graytext]  at (-0.2,4.9) {$v_6$};
  \node[vl,left]  at (-0.2,3.9) {$v_5$};
  \node[vl,left,text=graytext]  at (-1.85,2.55) {$v_3$};
  \node[vl,left,text=graytext]  at (1.45,2.55) {$v_8$};
  \node[vl,left,text=graytext]  at (-2.75,1.25) {$v_1$};
  \node[vl,left,text=graytext]  at (-0.65,1.25) {$v_4$};
  \node[vl,left,text=graytext]  at (-1.55,0.05) {$v_2$};
  \node[vl,left,text=graytext]  at (0.25,0.05) {$v_7$};
  \node[vl,right,text=graytext] at (2.95,2.55) {$v_{11}$};
  \node[vl,left,text=graytext]  at (1.45,1.25) {$v_{10}$};
  \node[vl,left,text=graytext]  at (1.45,0.05) {$v_9$};

\node[oldenc] at (0.5,5.3) {$[-1]$};
  %\node[encb] at (0.6,4.4) {$[[-1]]$};
  % \node[encb] at (-1.7,3.1) {$[[1]]$};
  % \node[encb] at (1.6,3.1) {$[[1]]$};
  \node[oldenc] at (-2.9,0.8) {$[1]$};
  % \node[encb] at (-0.2,1.8) {$[[-1],[1]]$};
  \node[oldenc] at (-1.7,-0.4) {$[-1]$};
  \node[oldenc] at (0.8,-0.4) {$[1]$};
  \node[oldenc] at (3.1,3.1) {$[1]$};
  %\node[encb] at (2.0,1.8) {$[[1]]$};
  \node[oldenc] at (2.0,-0.4) {$[1]$};

  % \node[oldenc] at (0.4,5.1) {$[-1]$};
  % \node[encb] at (0.6,4.4) {$[[-1]]$};
  %\node[encb] at (-2.3,3.1) {$[[1],[-3,-1,1]]$};
  %\node[encb] at (1.9,3.1) {$[[1],[2,1]]$};
  %\node[oldenc] at (-3.1,1.0) {$[1]$};
  \node[oldenc] at (-0.1,1.8) {$[-3,-1,1]$};
  %\node[oldenc] at (-1.7,-0.25) {$[-1]$};
  %\node[oldenc] at (0.7,-0.25) {$[1]$};
  %\node[oldenc] at (3.1,2.4) {$[1]$};
  \node[oldenc] at (2.0,1.8) {$[2,1]$};
  %\node[oldenc] at (2.0,-0.25) {$[1]$};

  %\node[oldenc] at (0.4,5.1) {$[-1]$};
  %\node[encb] at (0,4.4) {$[[-1],[5,-3,-1,1,1],[4,1,2,1]]$};
  \node[oldenc] at (-2.3,3.1) {$[5,-3,-1,1,1]$};
  \node[oldenc] at (1.9,3.1) {$[4,1,2,1]$};
  %\node[oldenc] at (-3.1,1.0) {$[1]$};
  %\node[oldenc] at (0.1,1.0) {$[-3,-1,1]$};
  %\node[oldenc] at (-1.7,-0.25) {$[-1]$};
  %\node[oldenc] at (0.7,-0.25) {$[1]$};
  %\node[oldenc] at (3.1,2.4) {$[1]$};
  %\node[oldenc] at (2.05,1.0) {$[2,1]$};
  %\node[oldenc] at (2.0,-0.25) {$[1]$};

  %\node[oldenc] at (0.4,5.1) {$[-1]$};
  \node[enc] at (0.3,4.4) {$[11,-1,4,1,2,1,5,-3,-1,1,1]=code(T)$};
  %\node[oldenc] at (-3.4,2.9) {$[5,-3,-1,1,1]$};
  %\node[oldenc] at (0.9,2.9) {$[4,1,2,1]$};
  %\node[oldenc] at (-3.1,1.0) {$[1]$};
  %\node[oldenc] at (0.1,1.0) {$[-3,-1,1]$};
  %\node[oldenc] at (-1.7,-0.25) {$[-1]$};
  %\node[oldenc] at (0.7,-0.25) {$[1]$};
  %\node[oldenc] at (3.1,2.4) {$[1]$};
  %\node[oldenc] at (2.05,1.0) {$[2,1]$};
  %\node[oldenc] at (2.0,-0.25) {$[1]$};
  
\node[
    font=\scriptsize,
    text=navy,
    align=center,
    text width=4cm
] at (0,-1)
    {(e) Final encoding};
\end{scope}

\end{tikzpicture}
\caption{The Encodings of Black and White Trees with 11 Nodes}
\label{Fig-04}
\end{figure}
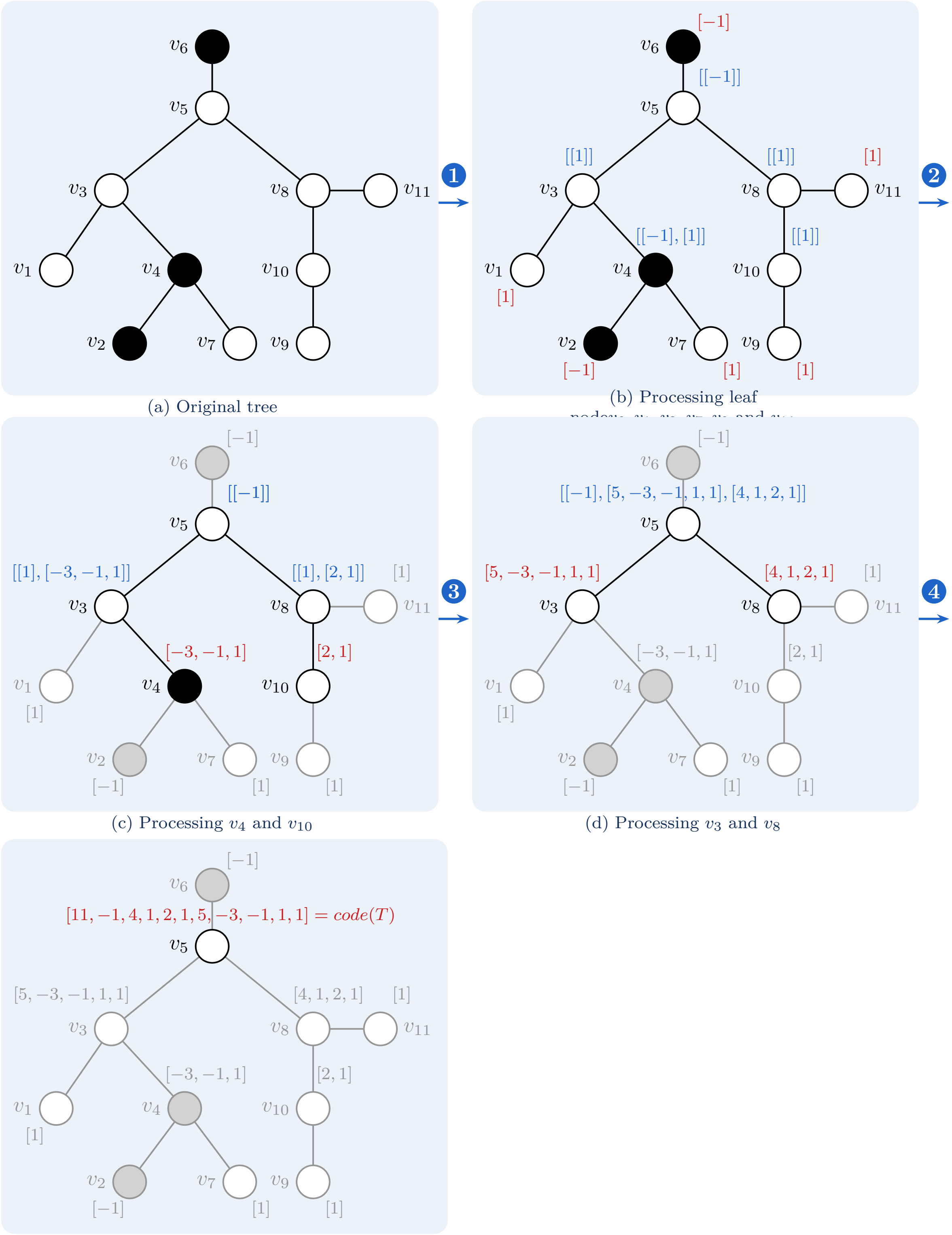

Following the approach used in the earlier definition of unrooted tree encoding, as outlined in \ref{defn210}, the encoding of unrooted trees with black and white nodes can be naturally extended. Below, we present an illustrative example.
\begin{example}
Consider the unrooted tree with black and white nodes shown in  \ref{Fig-04}. The encoding procedure initiates from the leaf nodes $v_1, v_2, v_6, v_7, v_9,$ and $v_{11}$, and proceeds as follows:
\end{example}
\begin{itemize}
    \item{Initialization:} The white leaf nodes $\{v_1, v_7, v_9, v_{11}\}$ are assigned the encoding $[1]$, while the black leaf nodes $\{v_2, v_6\}$ are encoded as $[-1]$. These values are propagated to their respective parent nodes $\{v_3, v_4, v_5, v_{10}, v_8\}$. Notably, since $v_2$ and $v_7$ share a common parent $v_4$, the parent node initially stores the aggregate encoding $[[-1], [1]]$.
    \item{Iteration 1:} Upon pruning the initial leaf nodes, a reduced subtree emerges where $v_4$ and $v_{10}$ become the new leaves.  For the black node $v_4$, we apply Definition 4 to derive the encoding $[-3, -1, 1]$. Here, the negative sign denotes the node’s black color, and the magnitude $3$ indicates the total number of nodes in the subtree rooted at $v_4$. The sequence $[-1, 1]$ represents the sorted, serialized encodings of its children. Similarly, the white node $v_{10}$ is encoded as $[2, 1]$. These results are passed to parents $v_3$ and $v_8$, which concatenate these new values with previously stored encodings (e.g., $[[1]]$). Consequently, the temporary encodings for $v_3$ and $v_8$ are updated to $[[1], [-3, -1, 1]]$ and $[[1], [2, 1]]$, respectively.
    \item{Iteration 2:}  Following the removal of the processed nodes, $v_3$ and $v_8$ become the leaf nodes of the remaining subtree. Applying Definition 4, we obtain the encodings $[5, -3, -1, 1, 1]$ for $v_3$ and $[4, 1, 2, 1]$ for $v_8$. These are transmitted to their common parent, $v_5$, which then stores the cumulative temporary encoding $[[-1], [5, -3, -1, 1, 1], [4, 1, 2, 1]]$.
    \item{Termination:} In the final stage, only node $v_5$ remains, representing both the leaf of the current subtree and the center node of the original tree $T$. By sorting the encodings of its children in ascending order and concatenating them, we derive the final representation of the original unrooted tree: $$\text{code}(T) = [11, -1, 4, 1, 2, 1, 5, -3, -1, 1, 1]$$
\end{itemize}

\section{Algorithm and Main Results}
\subsection{Algorithm Steps}
\begin{algorithm}[H]
\caption{Encoding a Tree, $\operatorname{code}(T)$}
\label{alg:code-tree}
\small                                   % <-- smaller body font
\SetAlgoSkip{smallskip}                  % <-- tighter vertical skips
\SetAlgoInsideSkip{smallskip}
\SetAlgoVlined
\DontPrintSemicolon
\KwIn{Adjacency matrix of a tree $T$}
\KwOut{Encoding $R$ of $T$}

\nonl
\textnormal{\footnotesize$\triangleright$ Initialize encoding variables}\;

$C \gets \operatorname{vector}(n,0)$\;
$R \gets [\,]$\;
\nonl
\textnormal{\footnotesize$\triangleright$ Prune leaves to the tree center}

\While{$\operatorname{remain}(T)>2$}{ 

    \ForEach{$i\in\operatorname{getLeaves}(T)$}{

        $C[i]\gets\operatorname{genCode}(C[i])$
        \tcp*[r]{Encode the leaf node $v_i$}

        $p\gets\operatorname{parent}(i,T)$
        \tcp*[r]{Find the parent of $v_i$}

        \textnormal{\footnotesize$\triangleright$ Update the parent's code} \\
        \eIf{$C[p]=0$}{
            $C[p]\gets\{C[i]\}$\;
        }{
            $C[p]\gets\{\operatorname{op}(C[p]),C[i]\}$\;
            $C[p]\gets\operatorname{sort}(C[p])$\;
        }
    }

    $T\gets\operatorname{removeLeaves}(T)$
    \tcp*[r]{Remove processed leaves}
}
\nonl
\textnormal{\footnotesize$\triangleright$ Encode the center node(s)}

\ForEach{$j\in\operatorname{getNodes}(T)$}{

    $C[j]\gets\operatorname{genCode}(C[j])$
    \tcp*[r]{Encode the center node $v_j$}

    \eIf{$\operatorname{nops}(R)=0$}{
        $R\gets C[j]$\;
    }{
        $R\gets\{op(R),C[j]\}$\;

        $R\gets\operatorname{sort}(R)$
    }
}

\Return{$R$}\;
\end{algorithm}

\indent Any tree with $n$ nodes can be constructed by appending a single node to an existing tree with $n-1$ nodes. This recursive property allows for the systematic organization of topological classifications for trees containing both black and white nodes. The specific procedure is as follows:

\begin{itemize}
    \item Case $n=1$: We start with the simplest case, in which the tree consists of a single node. Since the node can be either black or white, there are two possible colored trees, with encodings ([-1]) and ([1]), respectively.
    \item Case $n=2$: Starting from the two colored trees obtained for (n=1), we construct candidate trees by attaching a new node, either black or white, to an existing node. Although there are four possible combinations, applying the encoding and deduplication via Algorithm $code(T)$ yields only two distinct cases: $[-1, 1]$ and $[1, 1]$. Notably, $[-1, -1]$ is omitted, since the problem constraints forbid connecting two black nodes.
    \item Case $n=3$: We extend the colored trees obtained for (n=2) by attaching a new node, either black or white, to an existing node of each tree. Applying Algorithm $code(T)$ to identify isomorphic trees and eliminate duplicates, we obtain four distinct colored trees up to isomorphism, with encodings ([3,1,1]), ([3,-1,1]), ([3,-1,-1]), and ([-3,1,1]).
\end{itemize}

   By repeating this iterative expansion and deduplication process, we can determine the total number of heterogeneous tree classifications for any specified value of $n$.

   Due to the recursive architecture of the algorithm and the intensive deduplication of encodings, the process involves significant computational complexity. To address this, we optimized the implementation using both Maple and Python to handle encoding calculations efficiently.
   
   Under current optimizations, the algorithm successfully computes encodings for trees with up to $n=15$ nodes when considering both black and white node variations. In the simplified case where node colors are not introduced, the computational reach extends to $n=21$. The following two tables summarize the experimental results and corresponding execution times of our Python implementation.

  % Table 1: Colored trees
\begin{table}[H]
\raggedright
\caption{Classification Results and Computational Time for Trees with Black and White Nodes}
\label{tab:colored-trees}

\small
\setlength{\tabcolsep}{4pt}
\renewcommand{\arraystretch}{1.15}

\begin{tabular}{@{}l*{8}{c}@{}}
\toprule
Number of Nodes ($n$)
    & 2 & 3 & 4 & 5 & 6 & 7 & 8 & 9 \\
\midrule
Number of Tree Types
    & 2 & 4 & 10 & 26 & 75 & 234 & 768 & 2647 \\
Running Time (mm:ss)
    & 00:00 & 00:00 & 00:00 & 00:00
    & 00:00 & 00:00 & 00:00 & 00:00 \\
\midrule
Number of Nodes ($n$)
    & 10 & 11 & 12 & 13 & 14 & 15 & & \\
\midrule
Number of Tree Types
    & 9466 & 34818 & 131149 & 503640
    & 1965552 & 7777081 & & \\
Running Time (mm:ss)
    & 00:01 & 00:05 & 00:23 & 01:43
    & 06:34 & 25:16 & & \\
\bottomrule
\end{tabular}
\end{table}

% Table 2: Uncolored trees
\begin{table}[H]
\raggedright
\caption{Classification Results and Running Time for Uncolored Trees \quad \quad}
\label{tab:uncolored-trees}

\small
\setlength{\tabcolsep}{4pt}
\renewcommand{\arraystretch}{1.15}

\begin{tabular}{@{}l*{8}{c}@{}}
\toprule
Number of Nodes ($n$)
    & 2 & 3 & 4 & 5 & 6 & 7 & 8 & 9 \\
\midrule
Number of Tree Types
    & 1 & 1 & 2 & 3 & 6 & 11 & 23 & 47 \\
Running Time (mm:ss)
    & 00:00 & 00:00 & 00:00 & 00:00
    & 00:00 & 00:00 & 00:00 & 00:00 \\
\midrule
Number of Nodes ($n$)
    & 10 & 11 & 12 & 13 & 14 & 15 & 16 & 17 \\
\midrule
Number of Tree Types
    & 106 & 235 & 551 & 1301 & 3159 & 7741 & 19320 & 48629 \\
Running Time (mm:ss)
    & 00:00 & 00:00 & 00:00 & 00:00
    & 00:01 & 00:02 & 00:05 & 00:13 \\
\midrule
Number of Nodes ($n$)
    & 18 & 19 & 20 & 21 & & & & \\
\midrule
Number of Tree Types
    & 123867 & 317955 & 823065 & 2144505 & & & & \\
Running Time (mm:ss)
    & 00:33 & 01:29 & 04:04 & 11:25 & & & & \\
\bottomrule
\end{tabular}
\end{table}

\section{Discussion}
With the help of the generalized AHU algorithm promoted by the present article, we have completed a algorithmic classification of general n-dimensional Lotka-Volterra stably dissipative systems. And specific classification results were provided when $n \leq 15$. At the same time, we also corrected the erroneous results in existing literature regarding $n$=6, and 7. When n is relatively small ($n<8$), the dynamic behavior of low dimensional systems can be divided into the following four categories:

(1)	If the LaSalle's invariant set for a system contains only a positive equilibrium, then the system is globally stable.[\cite{Goh1977}, \cite{Hsu1978}]

(2)	If all the diagonal elements $a_{ii}$ of the interaction matrix $A=(a_{ij})$ are zero, then the system has non-canonical Hamiltonian formulation. [\cite{Plank1995}, \cite{Zhao2010}]

(3)	If a positive equilibrium is not unique, then every orbit has a limit depends on its initial value.[\cite{Redheffer1982}, \cite{Zhao2010}]

(4)	If the LaSalle's invariant sets of the systems are composed of differential algebraic equations, the dynamical behavior of the systems are determined by the structures of the LaSalle's invariant sets whether they contain periodic orbits or not.[\cite{Lu1989}, \cite{Yang2019}]

We speculate that for the general n-dimensional systems, there should not be any situations outside of the above four categories.

\printbibliography

\end{document}